\documentclass[reqno,12pt]{amsart}
\usepackage{ragged2e}
\usepackage{appendix}
\usepackage{amssymb}
\usepackage{amsfonts}
\usepackage{amsmath}
\usepackage{mathrsfs}
\usepackage{diagbox}
\usepackage{empheq}
\usepackage{color,soul}
\usepackage{bbm}
\usepackage{amsthm, graphicx,empheq}
\usepackage{framed}
\usepackage{booktabs}
\usepackage[backref,colorlinks,linkcolor=red,anchorcolor=green,citecolor=blue]{hyperref}
\usepackage{varwidth} % 引入 varwidth 包
\usepackage{lineno}  % show line number for draft. disable for final version
\newtheorem{theorem}{Theorem}[section]
\newtheorem{definition}{Definition}[section]
\newtheorem{remark}{Remark}[section]
\newtheorem{proposition}{Proposition}[section]
\numberwithin{equation}{section}

\newcommand{\mres}{\mathbin{\vrule height 1.6ex depth 0pt width
0.13ex\vrule height 0.13ex depth 0pt width 1.3ex}}

\makeatletter
\@namedef{subjclassname@20xx}{\textup{xxxx} Mathematics Subject Classification}
\makeatother
\begin{document}

\title[Hypersonic similarity law]{Hypersonic similarity law for steady compressible Euler flows past slender bodies within the framework of Radon measure solutions}

\author{Shifan Kang}
\author{Bingsong Long}
\author{Hairong Yuan}


\address[S. Kang]{Center for Partial Differential Equations, School of Mathematical Sciences,
East China Normal University, Shanghai 200241, China}\email{\tt 15332329608@163.com}

\address[B. Long]{School of Mathematics and Statistics, Huanggang Normal University, Hubei 438000, China}\email{\tt longbingsong@hgnu.edu.cn}

\address[H. Yuan]{School of Mathematical Sciences,  Key Laboratory of Mathematics and Engineering Applications (Ministry of Education) \& Shanghai Key Laboratory of PMMP,  East China Normal University, Shanghai 200241, China}\email{\tt hryuan@math.ecnu.edu.cn}


\date{\today}

\begin{abstract}
\iffalse 
In this paper, we establish a mathematical theory on statement and validation of the hypersonic similarity law within the framework of Radon measure solutions of steady compressible Euler equations. We consider two scenarios: 
(1) two-dimensional steady non-isentropic compressible Euler flows past an infinitely long slender curved wedge; 
(2) three-dimensional steady non-isentropic compressible Euler flows past an infinitely long axisymmetric cone.
 It turns out that, for the hypersonic flow passing through a slender body with tiny slenderness $\tau$, if the parameter $K\doteq M_{\infty}\tau$ is fixed, by taking  $\tau \to 0$ (i.e., the Mach number of the upcoming flow  $M_{\infty} \to \infty$), the flow field structures (after scaling)  no longer depend on the body's shape and the Mach number $M_{\infty}$ independently, but only on  $K$ and adiabatic index $\gamma$ of the polytropic gas.  Mathematically, for non-isentropic Euler flows, we find a new system of hypersonic small-disturbance equations to describe steady compressible hypersonic flows passing slender bodies.  We demonstrate that  as  $ \tau \to0$, under suitable  non-dimensional scalings, the Radon measure solutions of the original problems of hypersonic flow passing bodies converge to those  of corresponding hypersonic small-disturbance problems.  The explicit forms of the Radon measure solutions derived for the two scenarios effectively simplify the convergence analysis. 
 \fi

In this paper, we develop a mathematical theory for the statement and validation of the hypersonic similarity law within the framework of Radon measure solutions to the steady compressible Euler equations. We investigate two scenarios: 
(1) two-dimensional steady non-isentropic compressible Euler flows past infinitely long slender curved wedges, and 
(2) three-dimensional steady non-isentropic compressible Euler flows past infinitely long axisymmetric cones.
We find that for hypersonic flow over a slender body with a small slenderness parameter \(\tau\), if the parameter \(K \doteq M_{\infty} \tau\) is fixed, then as \(\tau \to 0\) (which corresponds to the Mach number of the incoming flow \(M_{\infty} \to \infty\)), the flow field structures, after scaling, become independent of the body's shape and the Mach number \(M_{\infty}\). Instead, they depend solely on \(K\) and the adiabatic exponent \(\gamma\) of the polytropic gas. 
Mathematically, we derive a new system of hypersonic small-disturbance equations to describe steady compressible hypersonic flows past slender bodies. We demonstrate that as \(\tau \to 0\), under suitable non-dimensional scalings, the Radon measure solutions of the original hypersonic flow problems converge to those of the corresponding hypersonic small-disturbance problems. The explicit forms of the Radon measure solutions obtained for the two scenarios facilitate the convergence analysis.

\end{abstract}

\allowbreak
\allowdisplaybreaks

\maketitle

\tableofcontents %disable for short paper

\section{Introduction}\label{sec1}

In engineering, usually a compressible flow with Mach number exceeding five is called hypersonic flow.  As the Mach number is increased to a high value, some specific flow characteristics become progressively more pronounced, such as the Mach number independence principle, thin shock layer, Newtonian-Busemann pressure law and the hypersonic similarity law (cf. \cite[Chapters 3 and 4]{Anderson}). The study of hypersonic flow has attracted considerable attention owing to its critical importance in aerodynamics \cite{Anderson,Andersonmodern,HypersonicBertin,hayeshft}. In this paper, we focus on a rigorous mathematical validation of the hypersonic similarity law within the framework of Radon measure solutions for steady non-isentropic compressible Euler flows over slender bodies, encompassing both two-dimensional curved wedges and three-dimensional axisymmetric cones. 

\iffalse 
One primary motivation for rigorous verification of the hypersonic similarity law is to establish simplified down-scaled models, recognized as the hypersonic small-disturbance problems, for the study of hypersonic flows around slender bodies. In addition, by considering Radon measure solutions, we can introduce the concentration boundary layers to overcome the difficulties encountered in the framework of Lebesgue measurable function solutions, such as when the Mach number of the flow is extremely large, some state variables (like density) go to infinity, and some (e.g. velocity) undergo dramatic changes.  Hence, problems related to hypersonic flow and hypersonic similarity law are of multi-scale in nature. Specific scale transformations obtained in this research can reveal physical phenomena that are difficult to detect at conventional scales, thereby enabling enhanced discernment of flow-field-quantity variations (see more in Remark \ref{remark multi-scaling}). 
\fi

A primary motivation for the rigorous verification of the hypersonic similarity law is to establish simplified down-scaled models, referred to as hypersonic small-disturbance problems, for studying hypersonic flows around slender bodies. Additionally, by considering Radon measure solutions, we can introduce concentration boundary layers to address challenges encountered within the framework of Lebesgue measurable function solutions. For instance, when the Mach number of the flow is extremely large, certain state variables (such as density) may approach infinity, while others (such as velocity) can exhibit dramatic changes. Consequently, problems related to hypersonic flow and the hypersonic similarity law are inherently multi-scale in nature. The specific scale transformations derived in this research can reveal physical phenomena that are difficult to detect at conventional scales, thereby enhancing our understanding of variations in flow-field quantities (see Remark \ref{remark multi-scaling} for further details).

In the following, we provide a more detailed background on hypersonic similarity laws, Radon measure solutions, and the contributions of this paper.

Let $M_{\infty}$ denote the Mach number of the oncoming flow, and $\tau$ the slenderness ratio of the body. Physically, the hypersonic similarity law means that for hypersonic flows with different Mach numbers $M_\infty$ past slender bodies with different slenderness ratios $\tau$, if the hypersonic similarity parameter $K\doteq M_{\infty}\tau$ and the adiabatic exponent $\gamma$ of the polytropic gas are fixed, then the flow structures are similar under appropriate scaling when $M_{\infty}$ is sufficiently large. The discovery of hypersonic similarity law dates back to 1946, originating from H. S. Tsien's work on steady irrotational flow \cite{Tsien}. The law was further extended to steady rotational flow  \cite{Hayes} and unsteady flow  \cite{Hamaker1953}. Interestingly, Hayes \cite{Hayes} found that the steady hypersonic flow over a slender body is equivalent to an unsteady flow in one less space dimension, which is also known as the {\em hypersonic equivalence principle} (see also \cite[Section 4.8]{Anderson}). In fact, both the hypersonic similarity law and the hypersonic equivalence principle appear in the development of hypersonic small-disturbance theory (see more in Hayes \cite{Hayes} and Van Dyke \cite{vandyke}). Recently, for two-dimensional straight or Lipschitz wedges, G.-Q. Chen et al. established a mathematical theory on the validation and the optimal convergence rate for the hypersonic similarity law in the context of BV (bounded variations) weak solutions to the steady compressible Euler equations  \cite{XW2024,XWstrate,XW2020,XW2023}. As for three-dimensional (axisymmetric) cones, the relevant mathematical result has not yet been found.

In the study of hypersonic flow passing an obstacle, the primary challenge lies in the appearance of special singularities within the flow field around the obstacle when the Mach number $M_{\infty}$ of oncoming flow is extremely high. These singularities make classical shock wave theory insufficient to describe the flow field. In fact, it has been observed that when $M_{\infty}$ approaches infinity, the shock surface becomes very close to the obstacle's surface (i.e., infinite-thin shock layer appears). Consequently, the mass, momentum and energy concentrate on the upwind boundary of the obstacle, causing the state quantities (such as mass density) to become infinite, which means the mass density is no longer a Lebesgue measurable function. Classical integrable weak solutions therefore fail to depict the hypersonic-limit flow field correctly, and one needs to introduce the notion of measure solutions (e.g., Radon measure solutions) to characterize the infinite-thin shock layers generated within hypersonic-limit flow. Then, it is natural to study the hypersonic similarity law within the framework of Radon measure solutions.

Mathematically, verifying the hypersonic similarity law under the framework of Radon measure solutions is to establish that the Radon measure solution of the original boundary value problem converges, after scaling, to the Radon measure solution of the corresponding scaled problem, one that is named {\em hypersonic small-disturbance problem}. It turns out that once proper definitions of Radon measure solutions, scaling, and convergence are given, we can demonstrate an explicit validation of the hypersonic similarity law. Notably, for non-isentropic Euler flows, we derived a new system of  hypersonic small-disturbance equations, and clarify that the hypersonic equivalence principle doe {\em not} hold for this system if we consider Radon measure solutions, which is quite different from the case in \cite[Section 4.4]{Anderson} and \cite{XW2024}; see Remarks \ref{rk h equivalence p1}--\ref{rk h equivalence p3} for details.

This work presents the first results on the hypersonic similarity law from the perspective of Radon measure solutions. One of the primary challenges is defining the similarity of flow structures for measure solutions that involve Dirac measures supported on the boundary of the body. Additionally, careful consideration must be given to the transformation of these Dirac measures under spatial transformations. These problems differ from the previous studies under the framework of weak solutions that are functions. In fact, G.-Q. Chen et al. established a comparison between the two entropy solutions of the compressible Euler equations with characteristic boundary conditions in \cite{XW2024,XWstrate,XW2020,XW2023} via the Glimm scheme or front tracking method. 
They established the stability of the global solution by perturbing the background solution and constructed the physical state near the boundary based on the assumption of small total variations of the physical boundary, which differs from the approach taken in this paper. In contrast, we examine the problem within the framework of Radon measure solutions, describing the flow field using measures that are absolutely continuous with respect to the Lebesgue measure in the flow region, along with weighted Dirac measures supported on the boundary of the body for concentration boundary layes. The advantage of this new approach is that it allows us to explicitly solve problems for a class of highly curved bodies, resulting in sharp outcomes that provide a solid foundation for further studies.

We remark that the problems of hypersonic flow past wedges or cones were successfully studied through Radon measure solutions in \cite{Jin2D,Q-W-Y2021,Yuan3DNewton,Yuan2020(strw)}. In these problems, the Newtonian-Busemann pressure law was justified by establishing proper Radon measure solutions. Note that the obstacles are fixed in these studies; namely, the slenderness ratio $\tau$ is fixed. Thus, the hypersonic similarity parameter $K$ is not held constant; specifically, $K\to \infty$ as $M_{\infty}\to \infty$. We also clarify in this paper that there are three types of hypersonic-limit flow. Further details are provided in Proposition \ref{prop 2D-Euler eq} and Remark \ref{Diff type for hl} below. Some related analysis on function solutions of the steady supersonic flow over wedges or axisymmetric cones can be found in \cite{GQchen3D,GQchenKuangjie2021,GQchenSlemMarshall,GQChenZYQ2006,CSX-XZP-YHC2002,supersonicflowand,HD-QAF2025,HD-ZYQ2019,LJ-WI-YHC2014} and the references therein.

The rest of the paper is structured into three sections. Section \ref{sec 2} mainly presents some preliminaries to the problem of supersonic flow past a two-dimensional curved wedge. We first rigorously formulate the corresponding mathematical problem, denoted as Problem A. After that, under the scaling (\ref{2D scaling})--(\ref{2D scaling E}) below, we present the corresponding hypersonic small-disturbance problem (denoted as Problem B), and state the main theorem on the case of two-dimensional wedges as Theorem \ref{main th for 2D}. Section \ref{sec 3} focuses on the proof of Theorem \ref{main th for 2D}. Specifically, we construct Radon measure solutions for both Problems A and B, and further validate the hypersonic similarity law through a comparative analysis of these solutions. Section \ref{sec 4} extends our investigation to the problem of hypersonic flow over a three-dimensional axisymmetric cone. We achieve a similar result by studying this problem in cylindrical coordinates. Appendix \ref{appendix a} provides an analysis of mathematical structure for the new system of  hypersonic small-disturbance equations (\ref{hyper sd eq})--(\ref{state eq for hsd}). Appendix \ref{appendixfactor} explains the influence of spatial scale transformations on the Radon measure solutions considered in this paper.

\section{Mathematical formulation of supersonic flow passing two-dimensional curved wedges}\label{sec 2}

\subsection{The problem of supersonic flow passing bodies}\label{sec 2.1}
The two-dimensional steady non-isentropic compressible Euler system is of the form:
  \begin{equation}\label{2D-Euler eq}
         \left\{\begin{array}{l} 
    \partial _x(\rho u)+\partial _y(\rho v)=0,\\  
    \partial _x(\rho u^2+p)+\partial _y(\rho uv)=0,\\
    \partial _x(\rho uv)+\partial _y(\rho v^2+p)=0,\\
    \partial _x(\rho uE)+\partial _y(\rho vE)=0,\\
        \end{array}\right. 
    \end{equation}
where $\rho$, $E$, and $(u,v)$ are unknowns and represent the density of mass, total enthalpy per unit mass, and velocity of the flow, respectively; $p$ is the scalar pressure. The constitutive relation (equation of state) is given by
    \begin{equation}\label{state eq}
      p=\frac{\gamma -1}{\gamma }\rho\cdot\big(E-\frac{1}{2}(u^2+v^2)\big),
    \end{equation}
with  $\gamma>1$ the adiabatic exponent of a polytropic gas. Let $c\doteq \sqrt{\gamma p/\rho}$ be the local sound speed of the gas, and $M\doteq\sqrt{u^2+v^2}/c$ the Mach number of the flow. For $M>1$, the flow is called supersonic, and it is well-known that (\ref{2D-Euler eq})--(\ref{state eq}) is then a hyperbolic system of conservation laws. 

 %-----------------------------------fig-----------
\begin{figure}[htb]
\centering
\includegraphics[scale=0.5]{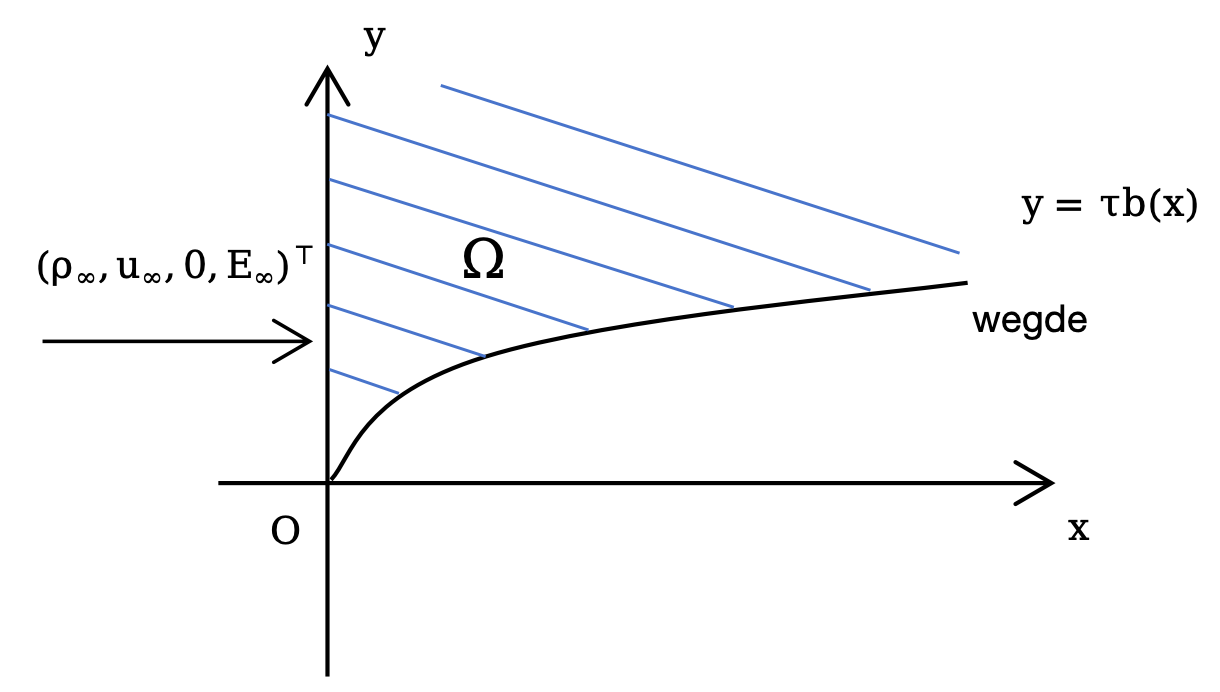}
\caption{Steady Euler flows over a two-dimensional slender wedge.}\label{fig2}
\end{figure}
 %-----------------------------------fig-----------
 
Now, let us describe our problem. On the $(x,y)$-plane, the slender wedge is composed of two boundaries given by the functions $y=\pm\tau b(x)$ for $x\ge0$, where $b(0)=0$ and $b'(x)>0$.  In this paper, the parameter $\tau$ satisfying $0<\tau\ll1$ is called  {\em slenderness ratio}. Set
  \begin{align*}
	\Omega\doteq\left \{ (x,y):~x \ge0,~y>\tau b(x)\right \},\quad \Gamma\doteq\{(x,y):~x\geq 0,~y=\tau b(x)\}.
  \end{align*}
By symmetry of the wedge, it is sufficient to study our problem in the domain $\Omega$ bounded by $\Gamma$ and the positive $y$-axis (see Figure \ref{fig2}). The flow satisfies the following slip condition
   \begin{equation}
        \label{slip boundary}
        (u,v)\cdot(n_1,n_2)=0\quad \text{on $\Gamma$},
    \end{equation} 
where $\textbf{n}\doteq (n_1,n_2)=(-\tau b'(x),1)/\sqrt{1+\tau^2b'(x)^2}$ is the unit normal vector on $\Gamma$ pointing to $\Omega$.

The oncoming flow of uniform state
   \begin{equation}
        \label{oncoming flow}
       U_\infty=(\rho_{\infty},u_{\infty},0,E_\infty)^{\top}
    \end{equation} 
is assumed to be supersonic, past the slender wedge without attack angle (namely $v_\infty=0$). For a fixed similarity parameter $K=M_{\infty}\tau$, there is 
\begin{equation*}
K^2=M_{\infty}^2\tau^2=\frac{u_{\infty}^2\tau^2}{\gamma\frac{p_{\infty}}{\rho_{\infty}}}=\frac {u_{\infty}^2\tau^2}{(\gamma-1)(E_{\infty}-\frac{1}{2}u_{\infty}^2)},
\end{equation*}
which implies
\begin{equation}
    \label{E p for onflow}
    \begin{aligned}
    E_{\infty} &= \frac{1}{2}u_{\infty}^2 + \frac{u_{\infty}^2 \tau^2}{(\gamma - 1)K^2} 
                = \frac{1}{2}u_{\infty}^2 \left(1 + \frac{2\tau^2}{(\gamma - 1)K^2}\right), \\
    p_{\infty} &= \frac{\rho_{\infty} u_{\infty}^2 \tau^2}{\gamma K^2}.
    \end{aligned}
\end{equation}
Therefore, our problem can be formulated mathematically as the following:
%%%%%%%%%%%%%%%%%%%%%%%%%%%%%%%%%%%%
\vspace{\baselineskip} % 段前间距    
\begin{center}    
\fbox{\begin{varwidth}{0.9\linewidth} % 设置宽度为行宽的80%        
\centering \textbf{Problem A}: For the oncoming flow $U_{\infty}$ given by (\ref{oncoming flow})--(\ref{E p for onflow}), find a solution to (\ref{2D-Euler eq})--(\ref{state eq}) in the domain $\Omega$ with the slip condition (\ref{slip boundary}). % 内容居中    
\end{varwidth}}    
\end{center}    
\vspace{\baselineskip} % 段后间距
%%%%%%%%%%%%%%%%%%%%%%%%%%%%%%%%%%%%%

We next derive the corresponding hypersonic small-disturbance problem to Problem A. Define the following non-dimensional independent and dependent variables:  
\begin{equation}\label{2D scaling}
    \begin{aligned}
     \bar{x}=x ,&\quad \bar{y}=\frac{y}{\tau} ,\\        \bar{u}=\frac{u-u_{\infty }}{u_{\infty}\tau^2},
    \quad \bar{v}=\frac{v}{u_{\infty}\tau},
   & \quad \bar{\rho}=\frac{\rho}{\rho_{\infty}},  \quad \bar{p}=\frac{p}{\gamma p_{\infty}M^2_{\infty}\tau^2},
    \end{aligned}
\end{equation}
see \cite[p.115]{Anderson} for the details of how these scalings were discovered from observations of Rankine-Hugoniot jump conditions. It follows from (\ref{state eq}) that a non-dimensional total enthalpy per unit mass can be defined as
\begin{equation}\label{2D scaling E}
     \bar{E}=\frac{2E-u_{\infty}^2}{u_{\infty}^2\tau^2}.
\end{equation}
Substituting (\ref{2D scaling})--(\ref{2D scaling E}) into (\ref{2D-Euler eq})--(\ref{state eq}), it yields
   \begin{equation}\label{2D Euler scaling}
        \left\{\begin{array}{l} 
        \partial _{\bar{x}}(\bar{\rho}( 1+\tau^2\bar{u})+\partial _{\bar{y}}(\bar{\rho} \bar{v})=0,\\  
        \partial _{\bar{x}}(\bar{\rho}\bar{u}( 1+\tau^2\bar{u})+\bar{p})+\partial _{\bar{y}}(\bar{\rho} \bar{u}\bar{v})=0,\\
        \partial _{\bar{x}}(\bar{\rho}\bar{v}( 1+\tau^2\bar{u}))+\partial _{\bar{y}}(\bar{\rho}\bar{v}^2+\bar{p})=0,\\
        \partial _{\bar{x}}(\bar{\rho}(1+\tau ^2\bar{u})\bar{E})+\partial _{\bar{y}}(\bar{\rho}\bar{v}
        \bar{E})=0,
        \end{array}\right. 
    \end{equation}
and
    \begin{equation}\label{state eq scaling}
        \bar{p}=\frac{\gamma-1}{2\gamma}\bar{\rho}(\bar{E}-2\bar{u}-\bar{v}^2-\tau^2\bar{u}^2 ).
    \end{equation}
Correspondingly, by (\ref{oncoming flow})--(\ref{2D scaling E}), the oncoming flow reads 
    \begin{equation}
        \label{oncoming flow scaling}
        \bar{U}_{\infty}=(1,0,0,\bar{E}_\infty)^{\top},
    \end{equation}
with
   \begin{equation}\label{E p for onflow scaling}
    \begin{aligned}
         \bar{E}_\infty&=\frac{2E_{\infty}-u_{\infty}^2}{u_{\infty}^2\tau^2}=\frac{2 }{(\gamma-1)K^2},\\
         \bar{p}_\infty&=\frac{p_{\infty}}{\gamma p_{\infty}M^2_{\infty}\tau^2}=\frac{1}{\gamma K^2}.
    \end{aligned}
    \end{equation}

On the $(\bar{x},\bar{y})$-plane, the domain occupied by the flow and its boundary are given respectively by
\begin{equation}\label{domain scaling}
    \Omega'=\left \{ (\bar{x},\bar{y}):~\bar{x} \ge0,~\bar{y}> b(\bar{x})\right \},\quad \Gamma' =\left \{ (\bar{x},\bar{y}):~\bar{x}\geq0,~\bar{y}= b(\bar{x}) \right \}.
\end{equation}
The boundary condition (\ref{slip boundary}) becomes
    \begin{equation}\label{slip boundary scaling}
       (1+\tau^2\bar{u},\bar{v})\cdot (\bar{n}_1,\bar{n}_2)=0\quad \text{on $\Gamma'$},
    \end{equation}
where $\bar{\textbf{n}}\doteq (\bar{n}_1,\bar{n}_2)=(- b'(\bar{x}),1)/\sqrt{1+b'(\bar{x})^2}$ is the unit normal vector on $\Gamma'$ pointing to $\Omega'$.
%Then, Problem A can be reformulated as
%\vspace{\baselineskip} % 段前间距    
%\begin{center}    
%\fbox{\begin{varwidth}{1\linewidth} % 设置宽度为行宽的80%        
%\centering \textbf{Problem A'}: For the oncoming flow $\bar{U}_{\infty}$ given above, find a solution to (\ref{2D Euler scaling}) in the domain $\Omega'$ with the slip boundary condition (\ref{slip boundary scaling}). % 内容居中    
%\end{varwidth}}    
%\end{center}    
%\vspace{\baselineskip} % 段后间距

\begin{figure}[htb]
\centering
\includegraphics[scale=0.5]{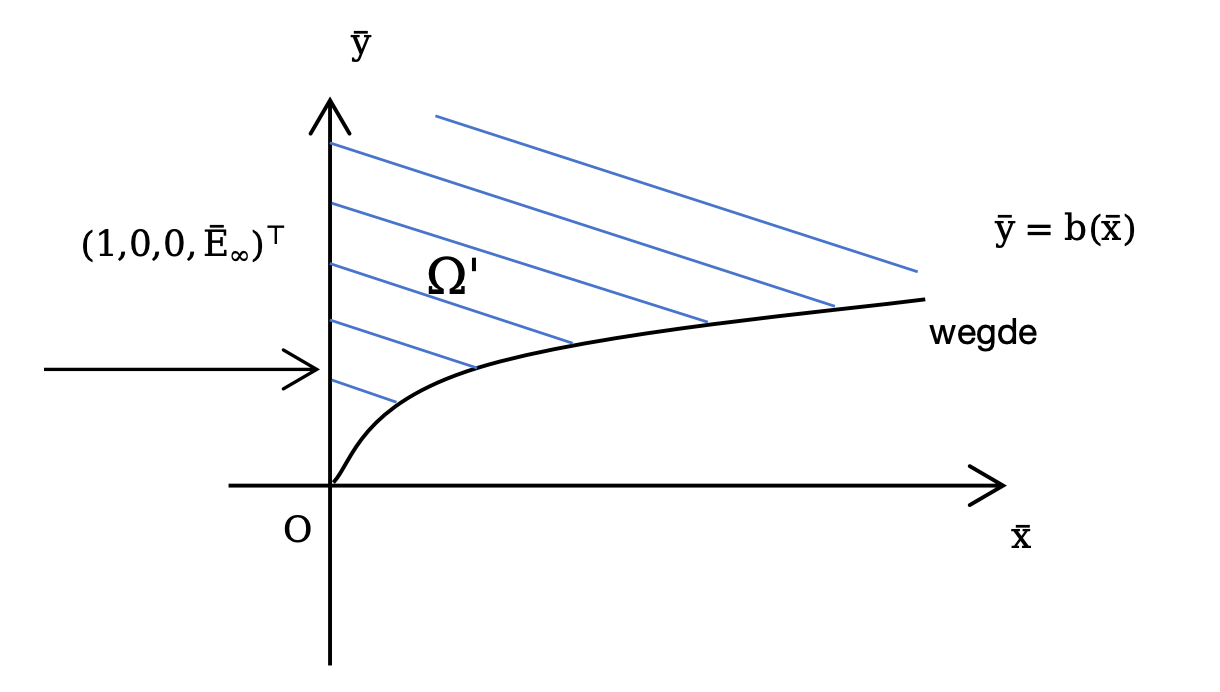}
\caption{{The corresponding two-dimensional hypersonic small-disturbance problem.}}\label{fig4}
\end{figure}

By the definition of $K$, the following is trivial. 

\begin{proposition}\label{prop 2D-Euler eq}
    For a fixed hypersonic similarity parameter~$K$, one has
\begin{equation*}
            \frac{1}{ M_{\infty }^2}=\frac{\tau^2}{K^2}.   
\end{equation*}
So the hypersonic limit $M_{\infty}\to \infty$ is equivalent to the vanishing slenderness ratio limit $\tau\to 0$.  
\end{proposition}
 
\begin{remark}\label{Diff type for hl}\rm
   In this paper, the hypersonic limit $\tau \to 0$ with the oncoming flow pressure $(\ref{E p for onflow scaling})_2$ differs from the classical hypersonic limit  which leads to the pressureless Euler flow as discovered in \cite{Yuan2020(strw)}. To be specific, for the case that $\tau>0$ is fixed, the hypersonic limit is corresponding to $K\to \infty$. From $(\ref{E p for onflow scaling})_1$, one has
   \begin{equation*}
       K^2=\frac{1}{(\gamma-1)\bar{E}_\infty}.
   \end{equation*}
Then there are other two cases of hypersonic limits:

\textit{Case} $1$. For fixed $\tau>0$ and (scaled) total enthalpy $\bar{E}_\infty>0$, the hypersonic limit $M_{\infty}\to \infty$ is equivalent to the limit $\gamma \to 1$, which corresponds to the vanishing pressure limit $($see, for example, \cite{QAF-YHR2020,Yuan2020(strw)}$)$;

\textit{Case} $2$. For fixed $\tau>0$ and $\gamma>1$, the hypersonic limit $M_{\infty}\to \infty$ is equivalent to the limit $\bar{E}_\infty \to 0$, which is not the classical vanishing pressure limit. In fact, there is no concentration in this case of hypersonic limit $($see, for example, \cite{HD-QAF2025,HD-ZYQ2019,Yuan2020(strw)}$)$. \hfill \qed
\end{remark}

For a fixed similarity parameter $K$, if the ratio $\tau\to 0$, then the terms involving $\tau^2$ in (\ref{2D Euler scaling})--(\ref{state eq scaling}) can be neglected intuitively; hence, we obtain the following {\em hypersonic small-disturbance equations}
\begin{equation}\label{hyper sd eq}
        \left\{\begin{array}{l} 
    \partial _{\bar{x}}\bar{\rho}+\partial _{\bar{y}}(\bar{\rho} \bar{v})=0,\\  
    \partial _{\bar{x}}(\bar{\rho}\bar{u}+\bar{p})+\partial _{\bar{y}}(\bar{\rho} \bar{u}\bar{v})=0,\\
     \partial _{\bar{x}}(\bar{\rho}\bar{v})+\partial _{\bar{y}}(\bar{\rho}\bar{v}^2+\bar{p})=0,\\
             \partial _{\bar{x}}(\bar{\rho}\bar{E})+\partial _{\bar{y}}(\bar{\rho}\bar{v}\bar{E})=0,
        \end{array}\right. 
    \end{equation}
with
    \begin{equation}\label{state eq for hsd}
        \bar{p}=\frac{\gamma-1}{2\gamma}\bar{\rho}(\bar{E}-2\bar{u}-\bar{v}^2).
    \end{equation}
Besides, the boundary condition (\ref{slip boundary scaling}) is reduced to
    \begin{equation}
        \label{slip condition for hsd}
        \bar{v}=b'(\bar{x})\quad \text{on $\Gamma'$}.
    \end{equation}

Then the hypersonic small-disturbance problem can be formulated as:
\vspace{\baselineskip} % 段前间距    
\begin{center}    
\fbox{\begin{varwidth}{0.9\linewidth} % 设置宽度为行宽的80%        
\centering \textbf{Problem B}: For the oncoming flow $\bar{U}_{\infty}$ given by (\ref{oncoming flow scaling})--(\ref{E p for onflow scaling}), find a solution to (\ref{hyper sd eq})--(\ref{state eq for hsd}) in the domain $\Omega'$ with the slip condition (\ref{slip condition for hsd}). % 内容居中    
\end{varwidth}}    
\end{center}    
\vspace{\baselineskip} % 段后间距
Note that Problem B  depends only on the parameter $K$ and  the adiabatic exponent $\gamma>1$.

\begin{remark}\label{rk h equivalence p1}\rm 
For steady supersonic Euler flow past a slender cone, Van Dyke \cite[p.4]{vandyke} and Anderson \cite[p.177]{Anderson} derived the classical hypersonic small-disturbance equations. By reducing the velocity component along the  $z$-axis, they simplify these equations to the following system in two-spatial-dimension:
   \begin{equation}\label{hsd eq vandyke}
    \left\{\begin{aligned} 
&\partial _{\bar{x}}\bar{\rho}+\partial _{\bar{y}}(\bar{\rho} \bar{v})=0,\\
&\bar{\rho} \partial _{\bar{x}}\bar{u} + \bar{\rho} \bar{v} \partial _{\bar{y}}\bar{u}= - \partial _{\bar{x}}\bar{p},
\\
&\bar{\rho} \partial _{\bar{x}} \bar{v} + \bar{\rho} \bar{v} \partial _{\bar{y}}\bar{v}= -\partial _{\bar{y}}\bar{p},
\\
&\partial _{\bar{x}}\left( \frac{\bar{p}}{\bar{\rho}^\gamma} \right) + \bar{v} \partial _{\bar{y}}\left( \frac{\bar{p}}{\bar{\rho}^\gamma} \right) = 0.
   \end{aligned}\right. 
    \end{equation}    
 Note that the first, third and fourth equations of $(\ref{hsd eq vandyke})$ are closed with respect to the (scaled) mass density function \(\bar{\rho}\), the velocity component \(\bar{v}\) and scalar pressure $\bar{p}$. This implies that $\bar{u}$ is decoupled from the system. Consequently, we can first solve for \(\bar{\rho}\), \(\bar{v}\) and $\bar{p}$, and then substitute them into the second equation of $(\ref{hsd eq vandyke})$ to determine the velocity component \(\bar{u}\). If we treat the spatial variable \(\bar{x}\) as the temporal variable \( t \), then the first, third and fourth equations of $(\ref{hsd eq vandyke})$ formally coincide with the one-dimensional unsteady compressible Euler equations. Furthermore, by subjecting the initial data (\ref{oncoming flow scaling})--(\ref{E p for onflow scaling}) and the boundary condition (\ref{slip condition for hsd}), it is the classical piston problem. This suggests that the hypersonic flow around a slender body can be approximately reduced to a simpler piston problem for analysis. This property, known as the \textit{hypersonic equivalence principle}, represents a significant application of the hypersonic small-disturbance equations. It should be pointed out that the fourth equation is deduced from the energy equation,  which only holds for classical solutions or isentropic flows. \hfill \qed
\end{remark}

\begin{remark}\label{rk h equivalence p2}\em 
For supersonic non-isentropic Euler flows involving shock waves, the entropy of flow exhibits a discontinuity $($jump$)$ across a shock front; namely, entropy is no longer constant along the streamline. Therefore, the system \eqref{hsd eq vandyke} is inadequate for describing such flows. To validate the hypersonic similarity law within the framework of BV weak solutions, G.-Q. Chen et al. established the following hypersonic small-disturbance equations \cite[$(1.12)$]{XW2024}
\begin{equation}\label{hsd eq G.Q.Chen}
        \left\{\begin{aligned} 
    &\partial _{\bar{x}}\bar{\rho}+\partial _{\bar{y}}(\bar{\rho} \bar{v})=0,\\  
   & \partial _{\bar{x}}(\bar{\rho}\bar{u}+\bar{p})+\partial _{\bar{y}}(\bar{\rho} \bar{u}\bar{v})=0,\\
     &\partial _{\bar{x}}(\bar{\rho}\bar{v})+\partial _{\bar{y}}(\bar{\rho}\bar{v}^2+\bar{p})=0,\\
    &\partial _{\bar{x}}(\bar{\rho}(\bar{u} + \frac{1}{2} \bar{v}^2 + \frac{\gamma \bar{p}}{(\gamma - 1) \bar{\rho}}))+\partial _{\bar{y}}(\bar{\rho}\bar{v}(\bar{u} + \frac{1}{2} \bar{v}^2 + \frac{\gamma \bar{p}}{(\gamma - 1) \bar{\rho}}))=0.
        \end{aligned}\right. 
    \end{equation}
Notably, the velocity component $\bar{u}$ can be decoupled from the other unknown variables in this system as well, and the hypersonic equivalence principle still holds.
\hfill\qed
\end{remark}

\begin{remark}\label{rk h equivalence p3}\rm 
Within the framework of Radon measure solutions, we shall adopt $\bar{\rho},\bar{u},\bar{v},\bar{E}$ as the primitive unknown variables, while treating pressure $\bar{p}$ as a derived physical quantity, given by $(\ref{state eq for hsd})$. Thus, for the hypersonic small-disturbance equations $(\ref{hyper sd eq})$ with $(\ref{state eq for hsd})$, the velocity component $\bar{u}$ cannot be decoupled, distinguishing it from both $(\ref{hsd eq vandyke})$ and $(\ref{hsd eq G.Q.Chen})$. The choice between $\bar{E}$ (the scaled specific enthalpy) and $\bar{p}$ (the scaled pressure) as a primitive variable is equivalent in the framework of function solutions; however, to investigate Radon measure solutions, selecting $\bar{E}$ as a primitive variable is more intrinsic due to its nature as a conserved quantity. This highlights a critical principle in the study of continuum physics: \textit{the choice of primitive unknown variables must align with the physical phenomena under investigation.}
\hfill\qed
\end{remark}

We point out here that Eqs. (\ref{hyper sd eq})--(\ref{state eq for hsd}) possess four real eigenvalues, including a double eigenvalue, rendering it a non-strictly hyperbolic system. It features two linearly degenerate characteristic fields and two genuinely nonlinear characteristic fields. The calculation is shown in Appendix \ref{appendix a}. Besides, (\ref{slip condition for hsd}) is a characteristic boundary condition.

\subsection{Main theorem on 2-d wedges}\label{sec 2.2}
Mathematically, the hypersonic similarity law is that the solution of Problem A, after scaling, converges in a suitable sense to that of Problem B as $\tau\to 0$.

For Problem A, we have the following result. 

   \begin{theorem}
        \label{th 2.1}
        Assume that $b(\cdot)\in C^2([0,\infty))$  satisfies $b(0)=0$, and for $x\ge0$, 
            \begin{align}
                 b'(x)&>0,\label{eq219}\\
                \frac{(1+\tau^2b'(x)^2)^{3/2}}{\gamma K^2}+b''(x)H(x)&+ b'(x)^2\sqrt[]{1+\tau^2b'(x)^2}>0, \label{eq220}
            \end{align}
        where
        \begin{equation}\label{eq for H}
           H(x)=\int_0^x\frac{ b'(t)}{\sqrt[]{1+\tau^2b'(t)^2}}\mathrm{d}t.
        \end{equation}
        Then, Problem A admits a Radon measure solution given by $(\ref{solution for A})$ below. Besides, for $x\ge0$, 
        \begin{equation}\label{NW law}
  w_p(x)=p_{\infty}+\frac{\rho_{\infty}u_{\infty}^2\tau^2 b''(x)H(x)+\rho_{\infty}u_{\infty}^2\tau^2 b'(x)^2\sqrt[]{1+\tau^2 b'(x)^2}}{(1+\tau^2b'(x)^2)^{3/2}}.
        \end{equation}
    \end{theorem}

\begin{remark}
    \label{remark:2.3}\rm 
    The formula $(\ref{NW law})$ is the {\em Newtonian-Busemann pressure law} for a wedge. In fact, by substituting $\theta$ with $\arctan{\tau b'(x)}$ and transforming the integration variable from $y$ to $x$, while utilizing the relation $\mathrm{d}y=\tau b'(x)\mathrm{d}x$, the expression presented in \cite[p.67, (3.27)]{Anderson} is equivalent to $(\ref{NW law})$. Moreover, the term $w_p(x)(-\mathbf{n})$ is of great interest in aerodynsmics, representing the force (lift/drag) exerted by the airflow on the wedge. Apparently, by $\eqref{E p for onflow}_2$ and \eqref{eq220}, we have $w_p(x)>0$, which means every single point on $\Gamma$ is subjected to the force of the flow. \hfill\qed
\end{remark}

As for Problem B, we have
\begin{theorem}\label{th2.2}
    Suppose that $b(\bar{x})\in C^2([0,\infty))$ satisfies $b(0)=0$, and for $x\ge0$, 
    \begin{equation*}
     b'(\bar{x}) > 0, \quad
     \bar{p}_\infty + b'(\bar{x})^2 + b(\bar{x})b''(\bar{x}) > 0.
    \end{equation*}
    Then Problem B admits a Radon measure solution provided by $(\ref{solution for hsd})$ below. Moreover, the corresponding term of $w_p(x)$ is 
     \begin{equation*}
        \bar{w}_p(\bar{x}) = \bar{p}_\infty + b'(\bar{x})^2 + b(\bar{x})b''(\bar{x}), \quad \bar{x}\ge0.
     \end{equation*}
    \end{theorem}

Our main result of the hypersonic similarity law for steady compressible Euler flow passing a two-dimensional slender wedge is the following theorem. 

\begin{theorem}[Main theorem for 2-d wedge]\label{main th for 2D}
Under the assumptions of Theorems \ref{th 2.1} and \ref{th2.2}, let ${\bar{U}^{(\tau)}}=(\bar{\varrho}^{(\tau)},\bar{u}^{(\tau)},\bar{v}^{(\tau)},\bar{E}^{(\tau)})^{\top}$ denote the Radon measure solution $U=(\varrho,u,v,E)^{\top}$ of Problem A after the scaling $(\ref{2D scaling})$--$(\ref{2D scaling E})$. Set $\mu=\mathcal{L}^2\mres\Omega'+\delta_{\Gamma'}$. Then when $\tau\to 0$, one has $\mu$-a.e. on $\Omega'\cup\Gamma'$ that 
        \begin{align}\label{eq223}
       & \bar{u}^{(\tau)}=\bar{u}+o(1), ~\bar{v}^{(\tau)}=\bar{v}+o(1), ~\bar{E}^{(\tau)}=\bar{E}+o(1), \\ 
       & \frac{\mathrm{d}\bar{\varrho}^{(\tau)}}{\mathrm{d}\mu}=\tau\frac{\mathrm{d}\bar{\varrho}}{\mathrm{d}\mu}+o(\tau), \label{eq224new}
    \end{align}
    and for each $\bar{x}\ge0$, 
    \begin{align}\label{eq225newa}
    w_p^{(\tau)}(\bar{x})=\bar{w}_p(\bar{x})+o(1),   \quad (\tau\to0),
     \end{align}
where $\bar{U}=(\bar{\varrho}, \bar{u}, \bar{v}, \bar{E})^\top$ is the Radon measure solution of Problem B, and $w_p^{(\tau)}(\bar{x})=w_p(\bar{x})/(\gamma p_\infty K^2)$.
\end{theorem}

The definition of Radon measure solutions, as well as the scaled solution ${\bar{U}^{(\tau)}}$ and the convergence \eqref{eq223}-\eqref{eq225newa} (particularly on the definition of the measure $\mu$ and the Radon-Nikodym derivatives ${\mathrm{d}\bar{\varrho}^{(\tau)}}/{\mathrm{d}\mu}$ etc.), will all be specified in the following section.

\section{Proof of main theorem for 2-d wedge}\label{sec 3}
In this section, after defining Radon measure solutions for Problems A and B, we will prove Theorem \ref{main th for 2D} to validate the hypersonic similarity law for the case of two-dimensional slender wedges. 

\subsection{Radon measure solution to Problem A}\label{sec  3.1}
We begin by reviewing some basic facts about Radon measures on the Euclidean space. Let $\mathcal{B}$ be the Borel $\sigma$-algebra of $\mathbb{R}^2$, and $m$ a  Radon measure on $(\mathbb{R}^2,\mathcal{B})$. By Riesz representation theorem, $m$ can be regarded as a continuous linear functional on $C_c(\mathbb{R}^2)$, i.e.,    
    \begin{equation}
        %\label{2.1}
        \left \langle m,\phi  \right \rangle =\int_{\mathbb{R}^2}\phi(x,y)m(\mathrm{d}x\mathrm{d}y),\quad \forall\,\text{$\phi\in C_c(\mathbb{R}^2)$},
     \end{equation}
where $C_c(\mathbb{R}^2)$ is the set of continuous functions on $\mathbb{R}^2$ with compact supports. Aside from the Lebesgue measure, another important example of Radon measure is the weighted Dirac measure.
%----------------------Df. Dirac measure-----------
    \begin{definition}[{Weighted Dirac measure supported on a curve}]
        \label{df:2.1}
        Let $$L\doteq \{ (x(t),y(t))\in\mathbb{R}^2:~ t\in[0,T] \}$$ be a $C^1$ curve. The weighted Dirac measure $w_L\delta_L$ supported on $L$ with a weight $w_L(\cdot)\in L_{loc}^1([0,T])$ is defined as 
        \begin{equation}
            \label{2.2}
        \left \langle w_L\delta_L,\phi  \right \rangle 
        =\int_0^Tw_L(t)\phi(x(t),y(t))\sqrt[]{x'(t)^2+y'(t)^2}\,\mathrm{d}t,\quad\forall\,\phi\in C_c(\mathbb{R}^2). 
        \end{equation}
\end{definition}

Before proceeding, we introduce some notations. Denote by $\chi_A$ the indicator function of a set $A$ (i.e. $\chi_A(x,y)=1$ for $(x,y)\in A$, and $\chi_A(x,y)=0$ otherwise); by $m\mres\Omega$ the measure obtained by restricting a measure $m$ on a Borel set $\Omega$; by $\mathcal{L}^d$ the standard Lebesgue measure of the Euclidean space $\mathbb{R}^d~(d\in\mathbb{N}$).
Also recall the standard notation $\lambda \ll \mu$ that a measure $\lambda$ is absolutely continuous with respect to a measure $\mu$. By Radon-Nikodym theorem,  if $\lambda$ and $\mu$ are $\sigma$-finite, there exists a $\mu$-a.e. unique function $\mathrm{d}\lambda/\mathrm{d}\mu$ called Radon-Nikodym derivative, such that
    \begin{equation*}
      \int_A \mathrm{d}\lambda= \int_A (\mathrm{d}\lambda/\mathrm{d}\mu)\mathrm{d}\mu, \qquad\forall\,A\in\mathcal{B}.
    \end{equation*}

\begin{definition}
    \label{df:2.2}
For a fixed $\gamma >1$, let $\varrho, \wp$ be nonnegative Radon measures, and $m^i, n^i~(i=0,1,2,3)$ be (signed) Radon measures on $\overline{\Omega}$, and $w_p(x)\in L_{loc}^1([0,\infty))$ be a nonnegative function. We call $(\varrho, u, v, E )$ a Radon measure solution to Problem A, provided that the following hold: 

    $(i)$$[$linear relaxation$]$ for any $\phi\in C_c^1(\mathbb{R}^2)$, there are
    \begin{equation}\label{2.3}
        \left \langle m^0 ,\partial_x\phi \right \rangle +
    \left \langle n^0,\partial_y\phi \right \rangle +
    \int _0^{\infty }\rho_{\infty}u_{\infty}\phi(0,y)\mathrm{d}y=0,
    \end{equation}
    \begin{equation}\label{2.4}
        \left \langle m^1+\wp ,\partial_x\phi \right \rangle +
    \left \langle n^1,\partial_y\phi \right \rangle +
    \left \langle w_p(x)n_1\delta _{y=\tau b(x)},\phi  \right \rangle +
    \int _0^{\infty }(\rho_{\infty }u_{\infty }^2+p_{\infty })\phi(0,y)\mathrm{d}y=0,
    \end{equation}
    \begin{equation}\label{2.5}
        \left \langle m^2,\partial_x\phi \right \rangle +
    \left \langle n^2+\wp,\partial_y\phi \right \rangle +
    \left \langle w_p(x)n_2\delta _{y=\tau b(x)},\phi  \right \rangle =0,
    \end{equation}
    \begin{equation}
        \label{2.6}
        \left \langle m^3 ,\partial_x\phi \right \rangle +
    \left \langle n^3,\partial_y\phi \right \rangle +
    \int _0^{\infty }\rho_{\infty}u_{\infty}E_{\infty}\phi(0,y)\mathrm{d}y=0,        
    \end{equation}
where $\mathbf{n}=(n_1,n_2)$ is the unit  normal vector on  $\Gamma$ pointing into $\Omega$;

$(ii)$$[$nonlinear constraints$]$ it is required that  $\wp \ll \varrho,~(m^k,n^k)\ll \varrho~~(k=0,1,2,3)$, with the Radon-Nikodym derivatives
\begin{equation}
       \label{2.7}
       u=\frac{\mathrm{d}m^0}{\mathrm{d}\varrho} \quad
       \mbox{and} \quad  v=\frac{\mathrm{d}n^0}{\mathrm{d}\varrho}
   \end{equation}
   satisfying $\varrho$-a.e. that
   \begin{equation}
       \label{2.8}
       uu=\frac{\mathrm{d}m^1}{\mathrm{d}\varrho} 
       ,~~uv=\frac{\mathrm{d}n^1}{\mathrm{d}\varrho}
            =\frac{\mathrm{d}m^2}{\mathrm{d}\varrho},~~
        vv=\frac{\mathrm{d}n^2}{\mathrm{d}\varrho},  
   \end{equation}   
   and there is a $\varrho$-a.e. function $E$ so that
      \begin{equation}
       \label{2.9}
        Eu=\frac{\mathrm{d}m^3}{\mathrm{d}\varrho},\quad
        Ev=\frac{\mathrm{d}n^3}{\mathrm{d}\varrho};
   \end{equation}
  set $u,v,E$ be zero on null set of $\varrho$;

$(iii)$$[$state equation$]$ if $\varrho \ll \mathcal{L}^2$~and~$\wp\ll \mathcal{L}^2$ on a Borel set $A$, then $\mathcal{L}^2\mres A$-a.e. there are
    \begin{equation*}
    \rho=\frac{\mathrm{d}\varrho}{\mathrm{d}\mathcal{L}^2},\quad 
      p=\frac{\mathrm{d}\wp}{\mathrm{d}\mathcal{L}^2}=\frac{\gamma -1}{\gamma }\rho\cdot\big(E-\frac{1}{2}(u^2+v^2)\big).
    \end{equation*}
Furthermore, classical entropy condition is valid for discontinuities of functions $\rho,u,v,E$ in this case.
\end{definition}

\begin{remark}
    \label{remark2.1} \rm
    We treat the compressible Euler equations as a linear differential system governed by the mass measure $\varrho$, momentum measures $m^0, n^0$, momentum flux measures $m^1, n^1, m^2, n^2$, enthalpy flux measures $m^3, n^3$, and pressure measure $\wp$. By employing the Radon-Nikodym theorem,  we recover the Euler system through the nonlinear relationships between the Radon-Nikodym derivatives of these measures. Since the Radon-Nikodym derivatives are defined point-wise, the challenges typically associated with the multiplication of measures in cases of concentration of physical quantities are not relevant here. The motivation and correctness of such a definition can be found in a series of works by the third author, see for instance, \cite{Jin2D, MR4773814,Q-W-Y2021,Yuan1dpitson,QAF-YHR2020,Yuan3DNewton,Yuan2020(strw)}. \hfill\qed
\end{remark}

\begin{remark}
    \label{remark:2.2} \rm
    Any  integral weak (function) solution can be considered naturally as  a Radon measure solution. The consistency of the definition of Radon measure solutions can be found in \cite{MR4773814,Yuan2020(strw)}.\hfill \qed
\end{remark}

Now we turn to construct a Radon measure solution for Problem A. It is natural to employ weighted Dirac measures supported on the boundary $\Gamma$ to represent the infinite-thin shock layers (i.e., concentration boundary layers). Thus, we suppose a Radon measure solution of Problem A is determined by the following regular-singular decompositions of measures with respect to the Lebesgue measure $\mathcal{L}^2$:

    \begin{equation}
        \label{2.10}
           m^0\doteq \rho_{\infty}u_{\infty }\mathcal{L}^2\mres\Omega+w_m^0(x)\delta_{y=\tau b(x)},
        \quad n^0\doteq w_n^0(x)\delta _{y=\tau b(x)},
    \end{equation}
    \begin{equation}
        \label{2.11}
           m^1\doteq \rho_{\infty}u_{\infty }^2\mathcal{L}^2\mres\Omega+w_m^1(x)\delta_{y=\tau b(x)},
        \quad n^1\doteq w_n^1(x)\delta _{y=\tau b(x)},\quad \wp\doteq p_{\infty}\mathcal{L}^2\mres\Omega,
    \end{equation}
    \begin{equation}
        \label{2.12}
           m^2\doteq w_m^2(x)\delta_{y=\tau b(x)}=w_n^1(x)\delta _{y=\tau b(x)},
        \quad n^2\doteq w_n^2(x)\delta _{y=\tau b(x)},
    \end{equation}
    \begin{equation}
        \label{2.13}
           m^3\doteq \rho_{\infty}u_{\infty }E_{\infty}\mathcal{L}^2\mres\Omega+w_m^3(x)\delta_{y=\tau b(x)},
    \quad n^3\doteq w_n^3(x)\delta _{y=\tau b(x)},
    \end{equation}
where $w_m^i(x), w_n^i(x)~~(i=0,1,2,3)$ are $C^1$ functions to be determined for $x\ge0$.

Inserting (\ref{2.10}) into (\ref{2.3}) gives
\begin{equation}\label{eq for zlsh}
\begin{aligned}
&\int_{\Omega}\rho_{\infty}u_{\infty}\partial_x\phi(x,y)\mathrm{d}x\mathrm{d}y\\
&\quad+\int_0^{\infty}\sqrt{1+\tau^2b'(x)^2}~\big[w_m^0(x)\partial
    _x\phi(x,\tau b(x))-w_n^0(x)\partial_y\phi(x,\tau b(x))\big]\mathrm{d}x\\&\qquad +
    \int _0^{\infty }\rho_{\infty}u_{\infty}\phi(0,y)\mathrm{d}y=0.
\end{aligned}
\end{equation}
It follows from Gauss-Green theorem that
\begin{equation}\label{eq for zlsh1} 
\begin{aligned}\int_{\Omega}\rho_{\infty}u_{\infty}\partial_x\phi(x,y)\mathrm{d}x\mathrm{d}y
&=\int_{\partial\Omega}\rho_{\infty}u_{\infty}\phi(x,y)\mathrm{d}y
\\&= \int _0^{\infty }\rho_{\infty}u_{\infty}\tau b'(x)\phi(x,\tau b(x))\mathrm{d}x-\int _0^{\infty }\rho_{\infty}u_{\infty}\phi(0,y)\mathrm{d}y.
\end{aligned}
\end{equation}
Using (\ref{eq for zlsh})--(\ref{eq for zlsh1}), and noticing that  
$$\partial_x\phi(x,\tau b(x))=\frac{\mathrm{d}}{\mathrm{d}x}\phi(x,\tau b(x))-\tau b'(x) \partial_y\phi(x,\tau b(x)), $$
we have
\begin{align*}
   & \int _0^{\infty }\Big(\rho_{\infty}u_{\infty}\tau b'(x)-\frac{\mathrm{d}(\sqrt[]{1+\tau^2b'(x)^2}~w_m^0(x))}{\mathrm{d}x}\Big)\phi(x,\tau b(x))\mathrm{d}x\\
    &\quad + \int _0^{\infty }\sqrt{1+\tau^2b'(x)^2}~\partial_y\phi(x,\tau b(x))(w_n^0(x)-\tau b'(x)w_m^0(x))\mathrm{d}x\\
    &\qquad -\sqrt{1+\tau^2b'(0)^2}~w_m^0(0)\phi(0,0)=0.
\end{align*}
This implies
\begin{equation}
    \label{2.14}
    \begin{aligned}
      \frac{\mathrm{d}(\sqrt[]{1+\tau^2b'(x)^2}~w_m^0(x))}{\mathrm{d}x}&=\rho_{\infty}u_{\infty}\tau b'(x),\\
      w_n^0(x)=\tau b'(x)w_m^0(x),&\quad w_m^0(0)=0,
    \end{aligned}
\end{equation}
due to the arbitrariness of the test funciton $\phi$.

Similarly, substituting (\ref{2.11})--(\ref{2.13}) into (\ref{2.4})--(\ref{2.6}) yields
\begin{equation}
    \label{2.15}
    \begin{aligned}
        \frac{\mathrm{d}(\sqrt[]{1+\tau^2b'(x)^2}~w_m^1(x))}{\mathrm{d}x}=(\rho_{\infty }u_{\infty}^2&+p_{\infty})\tau b'(x)-\tau b'(x)w_p(x),\\
         w_n^1(x)=\tau b'(x)w_m^1(x),&\quad w_m^1(0)=0,
    \end{aligned}
\end{equation}
and 
\begin{equation}
    \label{2.16}
    \begin{aligned}
        \frac{\mathrm{d}(\sqrt[]{1+\tau^2b'(x)^2}~w_m^2(x))}{\mathrm{d}x}&=-p_{\infty}+w_p(x),\\
         w_n^2(x)=\tau b'(x)w_m^2(x),&\quad w_m^2(0)=0,
    \end{aligned}
\end{equation}
as well as 
\begin{equation}
    \label{2.17}
    \begin{aligned}
        \frac{\mathrm{d}(\sqrt[]{1+\tau^2b'(x)^2}~w_m^3(x))}{\mathrm{d}x}&=\rho_{\infty }u_{\infty}E_{\infty}\tau b'(x),\\
         w_n^3(x)=\tau b'(x)w_m^3(x),&\quad w_m^3(0)=0.
    \end{aligned}
\end{equation}

With the use of (\ref{2.14})--(\ref{2.17}), we calculate the weight functions $w_m^i(x), w_n^i(x)~(i=0,1,2,3)$ and the pressure coeffcient $w_p(x)$. From (\ref{2.14}), a careful computation shows that
    \begin{equation}
        \label{2.18}
        w_m^0(x)=\frac{\rho_{\infty}u_{\infty}\tau b(x)}{\sqrt{1+\tau^2b'(x)^2}},\quad 
        w_n^0(x)=\frac{\rho_{\infty}u_{\infty}\tau^2 b(x)b'(x)}{\sqrt{1+\tau^2b'(x)^2}}.
    \end{equation}
We also solve (\ref{2.17}) to obtain
   \begin{equation}
        \label{2.19}
        w_m^3(x)=\frac{\rho_{\infty}u_{\infty}E_{\infty}\tau b(x)}{\sqrt{1+\tau^2b'(x)^2}},\quad 
        w_n^3(x)=\frac{\rho_{\infty}u_{\infty}E_{\infty}\tau^2 b(x)b'(x)}{\sqrt{1+\tau^2b'(x)^2}}.
    \end{equation}
From (\ref{2.12}) and $(\ref{2.15})_2$, it follows that
    \begin{equation*}
        \label{2.20}
        w_m^2(x)=w_n^1(x)=\tau b'(x)w_m^1(x).
    \end{equation*}
Combining this with $(\ref{2.15})_1$--(\ref{2.16}), one has
    \begin{equation}
    \label{2.24}
    \begin{aligned}
     w_n^1(x)&=w_m^2(x)=\frac{\rho_{\infty}u_{\infty}^2\tau^2 b'(x)H(x)}{1+\tau^2b'(x)^2}, \\
     w_m^1(x)&=\frac{\rho_{\infty}u_{\infty}^2\tau H(x)}{1+\tau^2b'(x)^2},\\
    w_n^2(x)&=\frac{\rho_{\infty}u_{\infty}^2 \tau^3b'(x)^2H(x)}{1+\tau^2b'(x)^2},
    \end{aligned}
    \end{equation}
    and
    \begin{equation}\label{2.25}
            w_p(x)=p_{\infty }+\frac{\rho_{\infty}u_{\infty}^2 \tau^2 b''(x)H(x)+\rho_{\infty}u_{\infty}^2 \tau^2 b'(x)^2\sqrt{1+\tau^2b'(x)^2}}{(1+\tau^2b'(x)^2)^{\frac{3}{2}}},
    \end{equation}
where $H(x)$ was given by (\ref{eq for H}).

By utilizing \eqref{2.10}--\eqref{2.13} and the equalities \eqref{2.18}--\eqref{2.24}, according to \eqref{2.7}–\eqref{2.9}, we can derive that
    \begin{equation}\label{2.27}
    \begin{aligned}
     u|_\Gamma&=\frac{u_{\infty} H(x)}{ b(x)\sqrt{1+\tau^2b'(x)^2}}=\frac{u_{\infty} \int_{0}^{x}\frac{ b'(t)}{\sqrt{1+\tau^2b'(t)^2}}\mathrm{d}t}{ b(x)\sqrt{1+\tau^2b'(x)^2}},  \\
    v|_\Gamma&=\frac{u_{\infty}\tau b'(x) H(x)}{ b(x)\sqrt{1+\tau^2b'(x)^2}}=\frac{u_{\infty}\tau b'(x) \int_{0}^{x}\frac{ b'(t)}{\sqrt{1+\tau^2b'(t)^2}}\mathrm{d}t}{ b(x)\sqrt{1+\tau^2b'(x)^2}},
    \end{aligned}
    \end{equation}
and
    \begin{equation}\label{2.28}
                E|_\Gamma=E_{\infty}.
    \end{equation}
Moreover, by (\ref{2.7}), (\ref{2.18}) and (\ref{2.27}), the singular part of $\varrho$ is
    \begin{equation}\label{eq326new}
        \frac{\rho_{\infty}\tau b(x)^2}{\int_0^x\frac{ b'(t)}{\sqrt{1+\tau^2b'(t)^2}}\mathrm{d}t}\delta_{y=\tau b(x)}.
    \end{equation}

From the above discussion, we infer that there exists the following Radon measure solution to Problem A:
   \begin{equation}
        \label{solution for A}
        \begin{aligned}
         &u=u_{\infty}\chi_{\Omega}+\frac{u_{\infty} \int_{0}^{x}\frac{ b'(t)}{\sqrt{1+\tau^2b'(t)^2}}\mathrm{d}t}{ b(x)\sqrt{1+\tau^2b'(x)^2}} ~\chi_\Gamma,\quad 
        v=\frac{u_{\infty}\tau b'(x) \int_{0}^{x}\frac{ b'(t)}{\sqrt{1+\tau^2b'(t)^2}}\mathrm{d}t}{ b(x)\sqrt{1+\tau^2b'(x)^2}}~\chi_\Gamma,\\
        &\varrho=\rho_{\infty}\mathcal{L}^2\mres\Omega+\frac{\rho_{\infty}\tau b(x)^2}{\int_0^x\frac{ b'(t)}{\sqrt{1+\tau^2b'(t)^2}}\mathrm{d}t}\delta_{y=\tau b(x)},\quad E=E_{\infty}\chi_{\Omega}+E_{\infty} \chi_\Gamma.
        \end{aligned}
    \end{equation}
It is straightforward to extend these expressions to $x=0$ by L'Hospital rule.   Thus, we proved Theorem \ref{th 2.1}. 

\begin{remark}\label{remark multi-scaling} \rm
It follows from \eqref{eq326new} that the weight function of $\varrho$ supported on the boundary $\Gamma$  vanishes as the slenderness ratio $\tau \to 0$, which corresponds to the vanishing boundary limits. However, for non-dimensional Problem~B $($corresponding to the vanishing limits of $\tau$$)$, we can observe the mass concentrated on the boundary $\Gamma'$ $($as to be shown in Section $\ref{sec 3.2})$. This is because the hypersonic similarity law is a typical multi-scale problem, which reveals physical phenomena that cannot be observed in the original large-scale formulation. \hfill \qed 
\end{remark}

\subsection{Radon measure solution to Problem B}\label{sec 3.2}
Under the scaling (\ref{2D scaling}), the Lebesgue measure $\mathcal{L}^2$ on $\Omega$ can be expressed as
\begin{equation}\label{eq df}
\begin{aligned}
   \left \langle \mathcal{L}^2\mres\Omega,\phi  \right \rangle &=\int_{\Omega}\phi(x,y)\mathrm{d}x\mathrm{d}y,\\
   &=\int_{\Omega'}{\bar{\phi}^{\tau}(\bar{x}, \bar{y})}(\tau\mathrm{d}\bar{x}\mathrm{d}\bar{y}),
\end{aligned}           
\end{equation}
where 
\begin{align}\label{eq328new}
\bar{\phi}^{\tau}(\bar{x},\bar{y})\doteq\phi(\bar{x}, \tau\bar{y})
\end{align} 
denotes the test function $\phi\in C_c(\mathbb{R}^2)$ under the scaling (\ref{2D scaling}). Since the validation of hypersonic similarity law is a multi-scale problem, we need to expand the observation with respect to the parameter $\tau$. To this end, the following definition is introduced. 
\begin{definition}
    \label{df2.3}
    Let $\bar{\mathcal{L}}^{2}\mres\Omega'$ be the Lebesgue measure restricted on the domain $\Omega'$, defined by 
    \begin{equation}\label{3.1}
    \left\langle \bar{\mathcal{L}}^{2}\mres\Omega', {\bar{\phi}}\right\rangle \doteq \int_{\Omega'}  \bar{\phi}(\bar{x}, \bar{y})\, \mathrm{d} {\bar{x}}\mathrm{d} {\bar{y}}, \qquad \forall\,\bar{\phi} \in C_c(\mathbb{R}^2).
    \end{equation}
    Also, let $L= \{ (\bar{x}(t),\bar{y}(t)): t\in[0,T] \} $ be a $C^1$ curve. Then the Dirac measure $\bar{w}_L\bar{\delta}_L$ supported on $L\subset \overline{\Omega'}$ with a weight $\bar{w}_L(\cdot)\in L_{loc}^1([0,T])$ is defined as 
    \begin{equation}\label{3.2}
    \left\langle \bar{w}_{L} \bar{\delta}_{L}, \bar{\phi}\right\rangle \doteq \int_{0}^{T} \bar{w}_{L}(t) \bar{\phi}(\bar{x}(t), \bar{y}(t))\sqrt{\bar{x}^{\prime}(t)^{2}+\bar{y}^{\prime}(t)^{2}}\, \mathrm{d} t, \qquad \forall\, \bar{\phi} \in C_c(\mathbb{R}^2).
    \end{equation}
    \end{definition}
\begin{remark}\label{remark dod for measure}\rm 
   Compared with $(\ref{eq df})$, the measure $(\ref{3.1})$ is defined by neglecting the scaling factor $\tau$ introduced by the transformation $\bar{x}=x,~\bar{y}=y/\tau$.  We proceed $(\ref{3.2})$ based on the consideration of $\bar{w}_{L}$ being the density of mass per unit length on the $(\bar{x},\bar{y})$-plane.   \hfill\qed 
\end{remark}

Correspondingly, we can define Radon measure solutions to Problem B.
\begin{definition}
    \label{df:2.4}
    For a fixed $\gamma> 1$, let $\bar{\varrho}, \bar{\wp}$ be nonnegative Radon measures, and $\bar{m}^i, \bar{n}^i~(i=0,1,2,3)$ be (signed) Radon measures on $\overline{\Omega'}$, and $\bar{w}_p(\bar{x})\in L_{loc}^1([0,\infty))$ be a nonnegative function. We call $(\bar{\varrho} ,\bar{u},\bar{v},\bar{E})$ a Radon measure solution to Problem B, if the following hold:

    $(i)[$linear relaxation$]$ for any $\bar{\phi}\in C_c^1(\mathbb{R}^2)$, there are
    \begin{equation}\label{2.33}
     \left \langle \bar{m}^0 ,\partial_{\bar{x}}\bar{\phi} \right \rangle + \left \langle \bar{n}^0,\partial_{\bar{y}}\bar{\phi} \right \rangle + \int _0^{\infty }\bar{\phi}(0,\bar{y})\mathrm{d}\bar{y}=0,
    \end{equation}
    \begin{equation}\label{2.34}
        \left \langle \bar{m}^1+\bar{\wp},\partial_{\bar{x}}\bar{\phi} \right \rangle +
    \left \langle \bar{n}^1,\partial_{\bar{y}}\bar{\phi} \right \rangle +\left \langle  \bar{w}_p(\bar{x})\bar{n}_1\bar{\delta}_{\Gamma'},\bar{\phi} \right \rangle 
    +\int _0^{\infty }\bar{p}_\infty\bar{\phi}(0,\bar{y})\mathrm{d}\bar{y}=0,
    \end{equation}
    \begin{equation}\label{2.35}
               \left \langle \bar{m}^2,\partial_{\bar{x}}\bar{\phi} \right \rangle +
    \left \langle \bar{n}^2+\bar{\wp},\partial_{\bar{y}}\bar{\phi} \right \rangle +\left \langle \bar{w}_p(\bar{x})\bar{n}_2\bar{\delta}_{\Gamma'},\bar{\phi} \right \rangle =0,
    \end{equation}
    \begin{equation}
        \label{2.36}
     \left \langle \bar{m}^3,\partial_{\bar{x}}\bar{\phi} \right \rangle +
    \left \langle \bar{n}^3,\partial_{\bar{y}}\bar{\phi} \right \rangle +\int_0^{\infty} \bar{E}_\infty\bar{\phi}(0,\bar{y})\mathrm{d}\bar{y}=0,   
    \end{equation}
where $\bar{\mathbf{n}}=(\bar{n}_1,\bar{n}_2)$ is the unit  normal vector on $\Gamma'$ pointing to $\Omega'$, and $\bar{m}^0=\bar{\varrho}$, $\bar{n}^0=\bar{m}^2$; 

$(ii)[$nonlinear constraints$]$ it is required that $\bar{\wp} \ll \bar{\varrho},~(\bar{m}^i,\bar{n}^i)\ll \bar{\varrho}~(i=0,1,2,3)$ with the Radon-Nikodym derivatives
\begin{equation}
       \label{2.37}    \bar{u}=\frac{\mathrm{d}\bar{m}^1}{\mathrm{d}\bar{\varrho}} \quad \text{and}\quad 
         \bar{v}=\frac{\mathrm{d}\bar{n}^0}{\mathrm{d}\bar{\varrho}}=\frac{\mathrm{d}\bar{m}^2}{\mathrm{d}\bar{\varrho}}
   \end{equation}
   satisfying $\varrho$-a.e. that
   \begin{equation}
       \label{2.38}    \bar{u}\bar{v}=\frac{\mathrm{d}\bar{n}^1}{\mathrm{d}\bar{\varrho}},\quad 
        \bar{v}\bar{v}=\frac{\mathrm{d}\bar{n}^2}{\mathrm{d}\bar{\varrho}}, 
   \end{equation}   
   and there is a $\bar{\varrho}$-a.e. function $\bar{E}$, so that
      \begin{equation}
       \label{2.39}     \bar{E}=\frac{\mathrm{d}\bar{m}^3}{\mathrm{d}\bar{\varrho}},\quad    \bar{E}\bar{v}=\frac{\mathrm{d}\bar{n}^3}{\mathrm{d}\bar{\varrho}};
   \end{equation}
  set $\bar{u},\bar{v},\bar{E}$ be zero on null set of $\bar{\varrho}$;

$(iii)[$state equation$]$ if $\bar{\varrho} \ll \bar{\mathcal{L}}^2$ and $\bar{\wp}\ll \bar{\mathcal{L}}^2$ on a Borel set, then $\bar{\mathcal{L}}^2$-a.e. there are
       \begin{equation*}
    \bar{\rho}=\frac{\mathrm{d}\bar{\varrho}}{\mathrm{d}\bar{\mathcal{L}}^2},\quad 
      \bar{p}=\frac{\mathrm{d}\bar{\wp}}{\mathrm{d}\bar{\mathcal{L}}^2}=\frac{\gamma-1}{2\gamma}\bar{\rho}(\bar{E}-2\bar{u}-\bar{v}^2).
    \end{equation*}
Moreover, classical entropy condition is valid for discontinuities of functions $\bar{\rho}, \bar{u},\bar{v},\bar{E}$ in this case.
\end{definition}

By the above definitions, we construct a Radon measure solution to Problem B according to the following decompositions:   
    \begin{equation}
        \label{3.10}       \bar{\varrho}=\bar{m}^0\doteq\bar{\mathcal{L}}^2\mres{\Omega'}+\bar{w}_{m}^0(\bar{x}) \bar{\delta}_{\bar{y}=b(\bar{x})},
        \quad \bar{n}^0\doteq\bar{w}_n^0(\bar{x}) \bar{\delta}_{\bar{y}=b(\bar{x})},
    \end{equation}
    \begin{equation}
        \label{3.11}
           \bar{m}^1\doteq \bar{w}_m^1(\bar{x}) \bar{\delta}_{\bar{y}=b(\bar{x})},
        \quad \bar{n}^1\doteq \bar{w}_n^1(\bar{x}) \bar{\delta}_{\bar{y}=b(\bar{x})},\quad \bar{\wp}\doteq \bar{p}_\infty\bar{\mathcal{L}}^2\mres{\Omega'},
    \end{equation}
    \begin{equation}
        \label{3.12}
          \bar{m}^2\doteq \bar{w}_m^2(\bar{x}) \bar{\delta}_{\bar{y}=b(\bar{x})}={\bar{w}_n^0(\bar{x})} \bar{\delta}_{\bar{y}=b(\bar{x})},\quad  \bar{n}^2\doteq\bar{w}_n^2(\bar{x}) \bar{\delta}_{\bar{y}=b(\bar{x})},
    \end{equation}
    \begin{equation}
        \label{3.13}      \bar{m}^3\doteq \bar{E}_\infty\bar{\mathcal{L}}^2\mres{\Omega'}+\bar{w}_{m}^3(\bar{x}) \bar{\delta}_{\bar{y}=b(\bar{x})},
        \quad \bar{n}^3\doteq \bar{w}_n^3 (\bar{x})\bar{\delta}_{\bar{y}=b(\bar{x})},
    \end{equation}
where the weights $\bar{w}_{m}^i(\cdot), \bar{w}_n^i(\cdot)~(i=0,1,2,3)$ are $C^1$ functions to be determined.

Proceeding as in Section \ref{sec  3.1}, we derive the ordinary differential equations obeyed by these weight functions:
    \begin{equation}
        \label{2.44}
      \left\{\begin{aligned}
             &\frac{\mathrm{d}(\sqrt[]{1+b'(\bar{x})^2}~\bar{w}_{m}^0(\bar{x}))}{\mathrm{d}\bar{x}}=b'(\bar{x}), \\
             &\frac{\mathrm{d}(\sqrt[]{1+b'(\bar{x})^2}~\bar{w}_m^1(\bar{x}))}{\mathrm{d}\bar{x}}=\bar{p}_\infty b'(\bar{x})- b'(\bar{x})\bar{w}_p(\bar{x}), \\
             &\frac{\mathrm{d}(\sqrt[]{1+b'(\bar{x})^2}~\bar{w}_m^2(\bar{x}))}{\mathrm{d}\bar{x}}=-\bar{p}_\infty +\bar{w}_p(\bar{x}), \\
             &\frac{\mathrm{d}(\sqrt[]{1+b'(\bar{x})^2}~\bar{w}_m^3(\bar{x}))}{\mathrm{d}\bar{x}}=\bar{E}_\infty b'(\bar{x}), \\
             \end{aligned}\right.
    \end{equation}
with the constraints 
    \begin{equation}
        \label{2.45}
        \begin{aligned}
       \bar{w}_n^i(\bar{x})= b'(\bar{x})\bar{w}_m^i(\bar{x}),\quad  \bar{w}_n^i(0)=0,\quad i=0,1,2,3.
         \end{aligned}
    \end{equation}
From (\ref{2.44})$_1$ and (\ref{2.45}), it follows that
    \begin{equation}
        \label{2.46}
        \bar{w}_m^0(\bar{x})=\frac{b(\bar{x})}{\sqrt[]{1+b'(\bar{x})^2} },\quad        \bar{w}_n^0(\bar{x})=\frac{b(\bar{x})b'(\bar{x})}{\sqrt[]{1+b'(\bar{x})^2} }, 
    \end{equation}
which means 
\begin{equation}\label{eqwm2}
    \bar{w}_m^2(\bar{x})=\frac{b(\bar{x})b'(\bar{x})}{\sqrt[]{1+b'(\bar{x})^2} },\quad  \bar{w}_n^2(\bar{x})=\frac{b(\bar{x})b'(\bar{x})^2}{\sqrt[]{1+b'(\bar{x})^2} }
\end{equation}
by the relation $\bar{w}_m^2(\bar{x})=\bar{w}_n^0(\bar{x})$ in (\ref{3.12}). Then, we can deduce from (\ref{2.44})$_2$, (\ref{2.44})$_3$, (\ref{2.45}) and (\ref{eqwm2}) that
     \begin{equation}
    \begin{aligned}
     \bar{w}_m^1(\bar{x})&=\frac{-b(\bar{x})b'(\bar{x})^2+\int_0^{\bar{x}}b(t)b'(t)b''(t)\mathrm{d}t }{\sqrt[]{1+b'(\bar{x})^2} },\\
     \bar{w}_n^1(\bar{x})&=\frac{-b(\bar{x})b'(\bar{x})^3+b'(\bar{x})\int_0^{\bar{x}}b(t)b'(t)b''(t)\mathrm{d}t}{\sqrt[]{1+b'(\bar{x})^2}},\\    \bar{w}_p(\bar{x})&=\bar{p}_\infty+b'(\bar{x})^2+b(\bar{x})b''(\bar{x})=\frac{1}{\gamma K^2}+b'(\bar{x})^2+b(\bar{x})b''(\bar{x}).
    \end{aligned}
    \end{equation}
Also, using (\ref{2.44})$_4$ and (\ref{2.45}), one has 
   \begin{equation}
        \label{2.47}
       \bar{w}_m^3(\bar{x})=\frac{\bar{E}_\infty b(\bar{x})}{\sqrt[]{1+b'(\bar{x})^2} }=\frac{\frac{2}{K^2(\gamma-1)} b(\bar{x})}{\sqrt[]{1+b'(\bar{x})^2} },\quad         \bar{w}_n^3(\bar{x})=\frac{\bar{E}_\infty b(\bar{x})b'(\bar{x})}{\sqrt[]{1+b'(\bar{x})^2} }=\frac{\frac{2}{K^2(\gamma-1)}b(\bar{x})b'(\bar{x})}{\sqrt[]{1+b'(\bar{x})^2} }  .
    \end{equation}
With (\ref{2.37})--(\ref{3.13}), the above relations imply that 
    \begin{equation}
    \begin{aligned}
     \bar{v}|_{\Gamma'}&=b'(\bar{x}), \\
     \bar{u}|_{\Gamma'}&=-(b'(\bar{x}))^2+\frac{\int_0^{\bar{x}}b(t)b'(t)b''(t)\mathrm{d}t}{b(\bar{x})},\\
     \bar{E}|_{\Gamma'}&=\bar{E}_\infty=\frac{2}{K^2(\gamma-1)}.
    \end{aligned}
    \end{equation}
    
Therefore, we obtain a Radon measure solution to Problem B as follows:
   \begin{equation}
        \label{solution for hsd}
        \begin{aligned}
        &\bar{u}=\Big(-b'(\bar{x})^2+\frac{\int_0^{\bar{x}}b(t)b'(t)b''(t)\mathrm{d}t}{b(\bar{x})}\Big) \chi_{\Gamma'},\quad \bar{v}=b'(\bar{x}) \chi_{\Gamma'} ,\\      &\bar{\varrho}=\bar{\mathcal{L}}^2\mres{\Omega'}+\frac{b(\bar{x})}{\sqrt[]{1+b'(\bar{x})^2} } \bar{\delta}_{\bar{y}=b(\bar{x})},\quad \bar{E}=\bar{E}_{\infty}\chi_{\Omega'}+\bar{E}_{\infty} \chi_{\Gamma'}.
        \end{aligned}
    \end{equation}
The extension of these expressions to $x=0$ is obvious. Theorem  \ref{th2.2} is proved.

\subsection{Comparison between Radon measure solutions for the two problems}\label{sec2.3}

In this subsection, we compare the Radon measure solutions obtained in Sections \ref{sec  3.1} and \ref{sec 3.2} to prove Theorem \ref{main th for 2D}, thus establishing the hypersonic similarity law. 

Recall the notation that ${\bar{U}^{(\tau)}(\bar{x})}=(\bar{\varrho}^{(\tau)},\bar{u}^{(\tau)},\bar{v}^{(\tau)},\bar{E}^{(\tau)})^{\top}$ denotes the Radon measure solution of Problem A under the scaling $(\ref{2D scaling})$--$(\ref{2D scaling E})$, and $\bar{U}(\bar{x})=(\bar{\varrho},\bar{u},\bar{v},\bar{E})^{\top}$ denotes the Radon measure solution of Problem B. The proof is divided into two steps.

\textit{Step} 1. {\em Compare the functions $\bar{u}^{(\tau)},\bar{v}^{(\tau)},\bar{E}^{(\tau)}$ and $\bar{u},\bar{v},\bar{E}$.} It is easy to justify that the solutions away from the boundary $\Gamma'$ are constants independent of $\tau$. Then, it remains to establish the convergence of the solutions on $\Gamma'$ when the slenderness ratio $\tau\to 0$. 

It follows from (\ref{2D scaling}) and (\ref{solution for A}) that
    \begin{equation*}
        \begin{aligned}
         \bar{u}^{(\tau)}|_{\Gamma'}=\frac{u|_\Gamma-u_{\infty}}{\tau^2u_{\infty}}= \frac{{\int_{0}^{\bar{x}}\frac{ b'(t)}{\sqrt{1+\tau^2b'(t)^2}}\mathrm{d}t}-b(\bar{x})\sqrt{1+\tau^2b'(\bar{x})^2}}{ \tau^2b(\bar{x})\sqrt{1+\tau^2b'(\bar{x})^2}}.  
          \end{aligned}
    \end{equation*}
Taking $\tau \to 0$, thanks to the  L'Hospital rule, we have
     \begin{equation}\label{2.57}
        \begin{aligned}
         \lim_{\tau\to 0} \bar{u}^{(\tau)}|_{\Gamma'}
          &=\lim_{\tau \to 0} \frac{{\int_{0}^{\bar{x}}\frac{ b'(t)}{\sqrt{1+\tau^2b'(t)^2}}\mathrm{d}t}-b(\bar{x})\sqrt{1+\tau^2b'(\bar{x})^2}}{ \tau^2b(\bar{x})}\\
          &=\lim_{\tau \to 0} \frac{{\int_{0}^{\bar{x}}\frac{-\tau  b'(t)^3}{(1+\tau^2b'(t)^2)^{3/2}}\mathrm{d}t}-\frac{\tau b(\bar{x})b'(\bar{x})^2}{\sqrt{1+\tau^2b'(\bar{x})^2}}}{ 2\tau b(\bar{x})}\\ 
             &= \lim_{\tau \to 0} \frac{{\int_{0}^{\bar{x}}\frac{-b'(t)^3}{(1+\tau^2b'(t)^2)^{3/2}}\mathrm{d}t}-\frac{b(\bar{x})b'(\bar{x})^2}{\sqrt{1+\tau^2b'(\bar{x})^2}}}{ 2b(\bar{x})}\\ 
            &= \frac{{\int_{0}^{\bar{x}}-  b'(t)^3\mathrm{d}t}-b(\bar{x})b'(\bar{x})^2}{2b(\bar{x})}.
          \end{aligned}
    \end{equation}
An integration-by-part yields 
        \begin{equation*} \begin{aligned}
        \int_{0}^{\bar{x}}- b'(t)^3\mathrm{d}t&=\int_{0}^{\bar{x}}- [(b(t)b'(t)^2)'-2b(t)b'(t)b''(t)]\mathrm{d}t\\
        &=-b(\bar{x})b'(\bar{x})^2+\int_{0}^{\bar{x}}2b(t)b'(t)b''(t)\mathrm{d}t.
     \end{aligned}    \end{equation*}
So (\ref{2.57}) can be written as
      \begin{equation}\label{2.58}
        \begin{aligned}
        \lim_{\tau\to 0} \bar{u}^{(\tau)}|_{\Gamma'}&=\frac{{-2b(\bar{x})b'(\bar{x})^2+\int_{0}^{\bar{x}}2b(t)b'(t)b''(t)\mathrm{d}t}}{2b(\bar{x})}\\
            &=-b'(\bar{x})^2+\frac{\int_{0}^{\bar{x}}b(t)b'(t)b''(t)\mathrm{d}t}{b(\bar{x})},
        \end{aligned}
        \end{equation}
which means
      \begin{equation*} 
        \lim_{\tau\to 0} \bar{u}^{(\tau)}|_{\Gamma'}=\bar{u}|_{\Gamma'}.
        \end{equation*}

Similarly, for the velocity along $\bar{y}$-axis, one has
\begin{equation*}
        \bar{v}^{(\tau)}|_{\Gamma'}=\frac{v|_\Gamma}{\tau u_{\infty}}
        =\frac{{b'(\bar{x})\int_{0}^{\bar{x}}\frac{ b'(t)}{\sqrt{1+\tau^2b'(t)^2}}\mathrm{d}t}}{ b(\bar{x})\sqrt{1+\tau^2b'(\bar{x})^2}}.
 \end{equation*}
As $\tau \to 0$,
\begin{align}\label{2.59}
        \lim_{\tau\to 0}\bar{v}^{(\tau)}|_{\Gamma'}
        = \frac{{b'(\bar{x})\int_{0}^{\bar{x}} b'(t)\mathrm{d}t}}{ b(\bar{x})}
        = b'(\bar{x})=\bar{v}|_{\Gamma'}.
 \end{align}
Besides, from the scaling (\ref{2D scaling E}), it follows that
\begin{align}
    \label{2.60}
    \lim_{\tau\to 0}\bar{E}^{(\tau)}|_{\Gamma'}=\lim_{\tau \to 0}\frac{2E|_\Gamma-u_{\infty}^2}{u_{\infty}^2\tau^2}=\bar{E}_\infty=\bar{E}|_{\Gamma'}.
\end{align}
The verification of $\lim_{\tau\to0}w_p^{(\tau)}(\bar{x})=\bar{w}_p(\bar{x})$ for $\bar{x}\ge0$ is also straightforward. 

\textit{Step} 2. {\em Compare the Radon measures $\bar{\varrho}^{(\tau)}$ and $\bar{\varrho}$.} To describe $\bar{\varrho}^{(\tau)}$, we recall that 
\begin{equation}\label{eq step2}
\begin{aligned}
    \left \langle \varrho,~\phi \right \rangle&=\left \langle\rho_{\infty}\mathcal{L}^2\mres{\Omega}+w_{\rho}(x)\delta_{y=\tau b(x)},~\phi\right \rangle\\
%    &=\left \langle\rho_{\infty}\mathcal{L}^2\mres{\Omega},~\phi\right \rangle+\left \langle w_{\rho}(x)\delta_{y=\tau b(x)},~\phi\right \rangle\\
    &=\int_{\Omega}\rho_{\infty}\phi~\mathrm{d}x\mathrm{d}y +\int_0^{\infty}w_{\rho}(x)\phi(x,\tau b(x))\sqrt{1+\tau^2b'(x)^2}~ \mathrm{d}x,
\end{aligned}
\end{equation}
where $w_{\rho}(x)$ denotes the weight function of the singular part of $\varrho$ (see (\ref{solution for A})), and $\phi \in C_c(\mathbb{R}^2)$ is a test function. Under the scaling $\bar{x}=x, \bar{y}=y/\tau$ and $\bar{\rho}=\rho/\rho_{\infty}$, the expression (\ref{eq step2}) can be reformulated as
\begin{equation}\label{2.61}
\begin{aligned}
 \left \langle\bar{\varrho}^{(\tau)},~\bar{\phi}^{(\tau)}\right \rangle&=\int_{\Omega'}\bar{\phi}^{(\tau)}~(\tau \mathrm{d}\bar{x}\mathrm{d}\bar{y})+\int_0^{\infty}\frac{w_{\rho}(\bar{x})}{\rho_{\infty}}\bar{\phi}^{(\tau)}(\bar{x}, b(\bar{x}))\sqrt{1+\tau^2b'(\bar{x})^2}~ \mathrm{d}\bar{x},\\
 &=\tau\Big(\int_{\Omega'}\bar{\phi}^{(\tau)}~\mathrm{d}\bar{x}\mathrm{d}\bar{y}+\int_0^{\infty}\frac{w_{\rho}(\bar{x})}{\tau \rho_{\infty}}\bar{\phi}^{(\tau)}(\bar{x}, b(\bar{x}))\sqrt{1+\tau^2b'(\bar{x})^2}~ \mathrm{d}\bar{x}\Big).
\end{aligned}
\end{equation}
Here $\bar{\phi}^{(\tau)}$ denotes $\phi$ under scaling (see \eqref{eq328new}), which is still arbitrary for any fixed $\tau>0$. Notice that the $w_{\rho}(\bar{x})$ is different from $\bar{w}_{\rho}(\bar{x})$ although they look similar, and the latter is $b(\bar{x})/\sqrt[]{1+b'(\bar{x})^2} $, the weight function of the singular part in  $\bar{\varrho}$ (see (\ref{solution for hsd})). 

The appearance of  the factor $\tau$ on the right hand side of \eqref{2.61} is explained in  Appendix \ref{appendixfactor}.  It is noted that \eqref{2.61} shall also be taken as a definition of the measure $\bar{\varrho}^{(\tau)}$.

We then consider $\bar{\varrho}$. For any test function $\bar{\phi}\in C_c(\mathbb{R}^2)$, 
\begin{equation}\label{2.62}
\begin{aligned}
    \left \langle\bar{\varrho},~\bar{\phi}\right \rangle&=\left \langle\bar{\mathcal{L}}^2\mres{\Omega'}+\bar{w}_{\rho}(\bar{x})\bar{\delta}_{\Gamma'},~\bar{\phi}\right \rangle\\
%    &=\left \langle\bar{\mathcal{L}}^2\mres{\Omega'},~\bar{\phi}\right \rangle+\left \langle\bar{w}_{\rho}(\bar{x})\bar{\delta}_{\Gamma'},~\bar{\phi}\right \rangle\\  
    &=\int_{\Omega'}\bar{\phi}~\mathrm{d}\bar{x}\mathrm{d}\bar{y} +\int_0^{\infty}\bar{w}_{\rho}(\bar{x})\bar{\phi}(\bar{x}, b(\bar{x}))\sqrt{1+b'(\bar{x})^2}~ \mathrm{d}\bar{x}.
\end{aligned}
\end{equation}
%\textcolor{red}{Of particular note is that we have expanded the observation on the scaling space $\bar{x}O\bar{y}$ with a parameter $\tau$ (see more in Appendix \ref{appendixfactor}), so it is natural to recover the difference with respect to $\tau$. } 
Obviously, by ignoring $\tau$ in (\ref{2.61}), the mass of gas $\bar{\varrho}^{(\tau)}$ per unit $\mathrm{d}\bar{x}\mathrm{d}\bar{y}$ on $\Omega'$ is  $1,$ 
which coincides with the mass of gas $\bar{\varrho}$ per unit $\mathrm{d}\bar{x}\mathrm{d}\bar{y}$ on $\Omega'$. From (\ref{2.61}), by neglecting the scaling factor $\tau$, the mass of gas $\bar{\varrho}^{(\tau)}$ per unit $\mathrm{d}\bar{x}$ on $\Gamma'$ is
\begin{equation*}
    \frac{w_{\rho}(\bar{x})}{\tau \rho_{\infty}}\sqrt{1+\tau^2b'(\bar{x})^2}.
\end{equation*}
To compare it  with the mass of gas $\bar{\varrho}$ per unit $\mathrm{d}\bar{x}$ on $\Gamma'$:
\begin{equation*}
    \bar{w}_{\rho}(\bar{x})\sqrt{1+b'(\bar{x})^2},
\end{equation*}
we use (\ref{solution for A}) and (\ref{solution for hsd}). It follows that 
\begin{equation}
    \lim_{\tau\to 0}\frac{\frac{w_{\rho}(\bar{x})}{\tau \rho_{\infty}}\sqrt{1+\tau^2b'(\bar{x})^2}}{ \bar{w}_{\rho}(\bar{x})\sqrt{1+b'(\bar{x})^2}}= \lim_{\tau\to 0}\frac{\frac{ b(\bar{x})^2\sqrt{1+\tau^2b'(\bar{x})^2}}{\int_0^{\bar{x}}\frac{ b'(t)}{\sqrt{1+\tau^2b'(t)^2}}\mathrm{d}t}}{b(\bar{x})}=1.
\end{equation}
This is what (\ref{eq224new}) means. The proof is completed.

\section{Hypersonic similarity law for three-dimensional Euler flows passing axisymmetric cones}\label{sec 4}

In this section, we establish analogous results on hypersonic similarity law for three-dimensional steady non-isentropic compressible supersonic Euler flows past  infinite-long axisymmetric cones. We first formulate the corresponding boundary value problem and hypersonic small-disturbance problem (denoted as Problem A$^{\prime}$ and Problem B$^{\prime}$, respectively), and then state the main theorem as Theorem \ref{main theorem 3D} in Section \ref{sec 4.1}. Owing to the axisymmetry of the cone, we study Problem A$^{\prime}$ and Problem B$^{\prime}$ in cylindrical coordinates, and treat these problems as a two-dimensional one, but with reference measures different from the Lebesgue measures. Subsequently, we construct the Radon measure solutions to Problem A$^{\prime}$ and Problem B$^{\prime}$ in cylindrical coordinates in Sections \ref{sec 4.3} and \ref{sec 4.4}. Finally, we prove the main theorem by comparing these Radon measure solutions constructed.

\subsection{Formulation of problems and main theorem}\label{sec 4.1}
\begin{figure}[htb]
\centering
\includegraphics[scale=0.5]{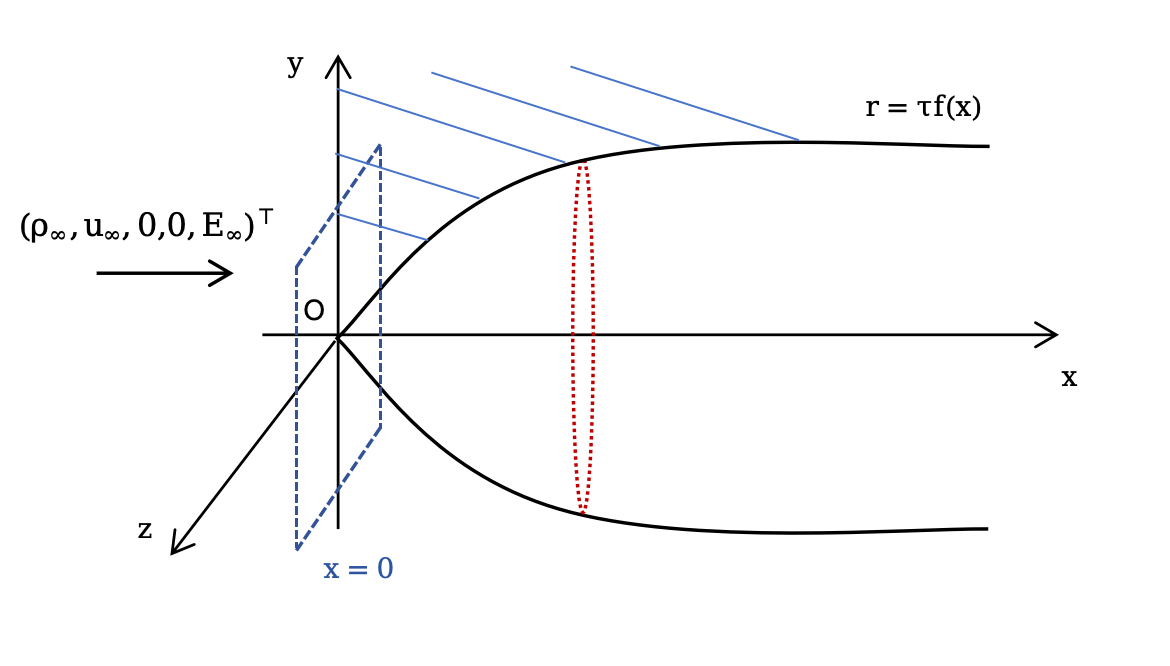}
\caption{Non-isentropic steady flow over an axisymmetric cone.}\label{fig5}
\end{figure}

We now formulate the mathematical problems for steady non-isentropic supersonic Euler flow passing a three-dimensional axisymmetric cone $$\left \{ (x,y,z)\in\mathbb{R}^3:~x\ge0,~0\le y^2+z^2\le\tau^2f(x)^2 \right \}$$ (see Figure \ref{fig5}), where $f(x)$ is a given $C^2$ function satisfying $f(0)=0$ and $f'(x)>0$ for $x\ge0$. Set 
\begin{equation}\label{domaineq}
    \Omega \doteq\left \{ (x,y,z)\in\mathbb{R}^3:~x\ge0,~y^2+z^2>\tau^2f(x)^2 \right \},
\end{equation}
and
\begin{equation}\label{boundaryeq}
    \Gamma\doteq\left \{ (x,y,z)\in\mathbb{R}^3:~x\ge0,~y^2+z^2= \tau^2f(x)^2 \right \}.
\end{equation}
Then, it suffices to study the problem in the domain $\Omega$.

The flow is governed by the following steady non-isentropic compressible Euler equations
    \begin{equation}\label{Euler_eq 3D}
        \left\{\begin{array}{l}
\partial_x(\rho u)+\partial_y(\rho v)+\partial_z(\rho w)=0,\\
\partial_x\left(\rho u^{2}+p\right)+\partial_y(\rho u v)+\partial_z(\rho u w)=0, \\
\partial_x(\rho u v)+\partial_y\left(\rho v^{2}+p\right)+\partial_z(\rho v w)=0, \\
\partial_x(\rho u w)+\partial_y(\rho v w)+\partial_z\left(\rho w^{2}+p\right)=0, \\
\partial_x(\rho u E)+\partial_y(\rho v E)+\partial_z(\rho w E)=0.
\end{array}\right.
    \end{equation}
Here $\rho$, $E$ are the density of mass and total enthalpy  per unit mass of the gas, respectively; and $u, v, w$ are the components of the velocity along $x,y$ and $z$-axis. The scalar pressure of the gas is given by   
\begin{equation}\label{state_eq 3D}
      p=\frac{\gamma -1}{\gamma }\rho\cdot\big(E-\frac{1}{2}(u^2+v^2+w^2)\big),
\end{equation}
where $\gamma>1$ is the adiabatic exponent of a polytropic gas. Besides, the flow satisfies
the slip condition
  \begin{equation}
        \label{slipcondition 3D}
        \tau uf(x)f'(x)=yv+zw \quad \text{on $\Gamma$}.
    \end{equation}
    
 %-----------------------------------fig5-----------
The uniform oncoming supersonic flow is given by
  \begin{equation}
        \label{oncoming flow 3D}
         U_{\infty}=(\rho_{\infty},u_{\infty},v_\infty\equiv0,w_\infty\equiv0,E_{\infty})^{\top}.
    \end{equation}
For fixed hypersonic similarity parameter $K$, a trivial calculation yields
\begin{equation}
    \label{oncoming E,p}
    \begin{aligned}
             E_{\infty}&=\frac{1}{2}u_{\infty}^2+\frac{u_{\infty }^2\tau^2}{(\gamma -1)K^2}=\frac{1}{2}u_{\infty}^2\Big(1+\frac{2\tau^2}{(\gamma -1)K^2}\Big),\\
    p_{\infty}&=\frac{\rho_{\infty}u_{\infty}^2\tau^2}{\gamma K^2}.
    \end{aligned}
\end{equation}
Therefore, we formulate the following boundary value problem.

%%%%%%%%%%%%%%%%%%%%%%%%%%%%%%%%%%%%
\vspace{\baselineskip}   
\begin{center}    
\fbox{\begin{varwidth}{0.9\linewidth} % 设置宽度为行宽的80%        
\centering \textbf{Problem A$^{\prime}$}: For the oncoming flow $U_{\infty}$ given by (\ref{oncoming flow 3D})--(\ref{oncoming E,p}), find a solution to (\ref{Euler_eq 3D})--(\ref{state_eq 3D}) in the domain $\Omega$ with the slip condition (\ref{slipcondition 3D}).    
\end{varwidth}}    
\end{center}    
\vspace{\baselineskip} % 段后间距
%%%%%%%%%%%%%%%%%%%%%%%%%%%%%%%%%%%%%

\begin{figure}[htb]
\centering
\includegraphics[scale=0.6]{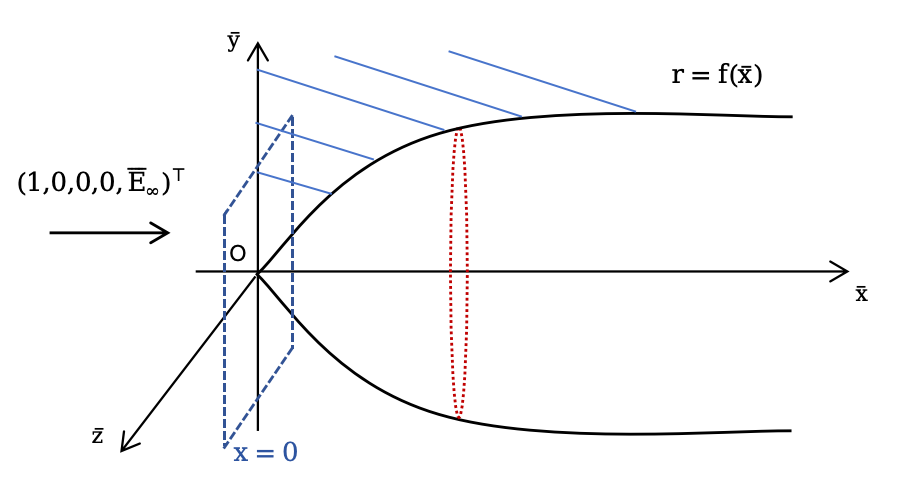}
\caption{Corresponding three-dimensional hypersonic small-disturbance problem.}\label{fig6}
\end{figure}

We next describe the corresponding hypersonic small-disturbance problem. Following \cite[p.115]{Anderson}, we introduce the scaling:
\begin{equation}
    \label{scaling 3D}
    \begin{aligned}
            \bar{x}=x ,\quad\bar{y}=\frac{y}{\tau},&\quad\bar{z}=\frac{z}{\tau},\quad\bar{\rho}=\frac{\rho}{\rho_{\infty}},\\
   \quad \bar{u}=\frac{u-u_{\infty }}{u_{\infty}\tau^2},
    \quad\bar{v}=\frac{v}{u_{\infty}\tau},
     &\quad\bar{w}=\frac{w}{u_{\infty}\tau},
     \quad \bar{p}=\frac{p}{\gamma p_{\infty}M^2_{\infty}\tau^2},
    \end{aligned}
\end{equation}
and by (\ref{state_eq 3D}), one has
\begin{equation}
    \label{scaling_E 3D}
    \bar{E}=\frac{2E-u_{\infty}^2}{\tau^2u_{\infty}^2}.
\end{equation}
Substituting (\ref{scaling 3D})--(\ref{scaling_E 3D}) into (\ref{Euler_eq 3D})--(\ref{state_eq 3D}), we obtain
\begin{equation}\label{3.8}
        \left\{\begin{array}{l} 
        \partial _{\bar{x}}(\bar{\rho}( 1+\tau^2\bar{u}))+\partial _{\bar{y}}(\bar{\rho} \bar{v})+\partial _{\bar{z}}(\bar{\rho} \bar{w})=0,\\ 
        \partial _{\bar{x}}(\bar{\rho}\bar{u}( 1+\tau^2\bar{u})+\bar{p})+\partial _{\bar{y}}(\bar{\rho} \bar{u}\bar{v})+\partial _{\bar{z}}(\bar{\rho} \bar{u}\bar{w})=0,\\
        \partial _{\bar{x}}(\bar{\rho}\bar{v}( 1+\tau^2\bar{u}))+\partial _{\bar{y}}(\bar{\rho}\bar{v}^2+\bar{p})+\partial _{\bar{z}}(\bar{\rho}\bar{v}\bar{w})=0,\\
        \partial _{\bar{x}}(\bar{\rho}\bar{w}( 1+\tau^2\bar{u}))+\partial _{\bar{y}}(\bar{\rho}\bar{v}\bar{w})+\partial _{\bar{z}}(\bar{\rho}\bar{w}^2+\bar{p})=0,\\
        \partial _{\bar{x}}(\bar{\rho}(1+\tau ^2\bar{u})\bar{E})+\partial _{\bar{y}}(\bar{\rho}\bar{v}
        \bar{E})+\partial _{\bar{z}}(\bar{\rho}\bar{w}
        \bar{E})=0,
        \end{array}\right. 
    \end{equation}
and
    \begin{equation}\label{3.9}
        \bar{p}=\frac{\gamma-1}{2\gamma}\bar{\rho}\cdot\big(\bar{E}-2\bar{u}-\bar{v}^2-\bar{w}^2-\tau^2\bar{u}^2 \big).
    \end{equation}
In the $(\bar{x},\bar{y},\bar{z})$-coordinates, the corresponding fluid domain and its boundary are given respectively 
by  (see Figure \ref{fig6})
\begin{align}\label{domainboundaryeq}
    &\Omega'=\left \{ (\bar{x},\bar{y},\bar{z}):~\bar{x} \ge0,~\bar{y}^2+\bar{z}^2>f(\bar{x})^2\right \},\\ 
    &{\Gamma'} =\{ (\bar{x},\bar{y},\bar{z}):~\bar{x} \ge0,~\bar{y}^2+\bar{z}^2= f(\bar{x})^2 \}.
\end{align}
The given oncoming flow (\ref{oncoming flow 3D})--(\ref{oncoming E,p}) becomes
    \begin{equation}
        \label{oncoming flow_hsd 3D}
        \bar{U}_{\infty}=(1,0,0,0,\bar{E}_\infty)^{\top}, \quad  \bar{E}_\infty=\frac{2}{(\gamma-1)K^2},\quad \bar{p}_\infty=\frac{1}{\gamma K^2}, 
    \end{equation}
and the slip condition (\ref{slipcondition 3D}) now reads  
    \begin{equation}\label{slip boundary_af 3D}     (1+\tau^2\bar{u})f(\bar{x})f'(\bar{x})=\bar{v}\bar{y}+\bar{w}\bar{z}.
    \end{equation}

 %-----------------------------------fig6-----------

From Proposition \ref{prop 2D-Euler eq}, as $M_{\infty}\to \infty$, the slenderness ratio $\tau\to 0$. Thus we intuitively neglect $\tau^2$ in (\ref{3.8})--(\ref{3.9}) to obtain the corresponding hypersonic small-disturbance equations 
    \begin{equation}\label{hsd govern eq 3D}
        \left\{\begin{array}{l} 
    \partial _{\bar{x}}\bar{\rho}+\partial _{\bar{y}}(\bar{\rho} \bar{v})+\partial _{\bar{z}}(\bar{\rho} \bar{w})=0,\\  
    \partial _{\bar{x}}(\bar{\rho}\bar{u}+\bar{p})+\partial _{\bar{y}}(\bar{\rho} \bar{u}\bar{v})+\partial _{\bar{z}}(\bar{\rho} \bar{u}\bar{w})=0,\\
     \partial _{\bar{x}}(\bar{\rho}\bar{v})+\partial _{\bar{y}}(\bar{\rho}\bar{v}^2+\bar{p})+\partial _{\bar{z}}(\bar{\rho} \bar{v}\bar{w})=0,\\
     \partial _{\bar{x}}(\bar{\rho}\bar{w})+\partial _{\bar{y}}(\bar{\rho} \bar{v}\bar{w})+\partial _{\bar{z}}(\bar{\rho}\bar{w}^2+\bar{p})=0,\\
             \partial _{\bar{x}}(\bar{\rho}\bar{E})+\partial _{\bar{y}}(\bar{\rho}\bar{v}\bar{E})+\partial _{\bar{z}}(\bar{\rho}\bar{w}\bar{E})=0,
        \end{array}\right. 
    \end{equation}
with
    \begin{equation}\label{hsd state eq 3D}
        \bar{p}=\frac{\gamma-1}{2\gamma}\bar{\rho}\cdot\big(\bar{E}-2\bar{u}-\bar{v}^2-\bar{w}^2\big).
    \end{equation}
Besides, the boundary condition (\ref{slip boundary_af 3D}) becomes
    \begin{equation}
        \label{slip boundary_hsd 3D}
        f(\bar{x})f'(\bar{x})=\bar{v}\bar{y}+\bar{w}\bar{z}.
    \end{equation}
Then we present the corresponding hypersonic small-disturbance problem as
 %%%%%%%%%%%%%%%%%%%%%%%%%%%%%%%%%%%%
\vspace{\baselineskip} % 段前间距    
\begin{center}    
\fbox{\begin{varwidth}{0.9\linewidth} % 设置宽度为行宽的80%        
\centering \textbf{Problem B$^{\prime}$}: For the oncoming flow $\bar{U}_{\infty}$ given by (\ref{oncoming flow_hsd 3D}), find a solution to (\ref{hsd govern eq 3D})--(\ref{hsd state eq 3D}) in the domain $\Omega'$ with the slip condition (\ref{slip boundary_hsd 3D}). % 内容居中    
\end{varwidth}}    
\end{center}    
\vspace{\baselineskip} % 段后间距
%%%%%%%%%%%%%%%%%%%%%%%%%%%%%%%%%%%%%

We will show that Problem A$^{\prime}$ and Problem B$^{\prime}$ admit axisymmetric Radon measure solutions provided that the boundary of the axisymmetric cone satisfies appropriate conditions (see Theorems \ref{th3.2} and \ref{th3.3} below). For steady supersonic Euler flows over an infinite-long axisymmetric cone, the main result can be stated as follows.

\begin{theorem}[Main theorem for 3-d axisymmetric cone]\label{main theorem 3D}\label{thm41}
Let $$\bar{U}^{(\tau)}=(\bar{\varrho}^{(\tau)},\bar{u}^{(\tau)},\bar{v}^{(\tau)},\bar{w}^{(\tau)},\bar{E}^{(\tau)})^{\top}$$ denote the axisymmetric Radon measure solution $U=(\varrho,u,v,w,E)^{\top}$ of $Problem~A^{\prime}$ under the scaling $(\ref{scaling 3D})$--$(\ref{scaling_E 3D})$. Set $\bar{\mu}=\mathcal{L}^3\mres\Omega'+\delta_{\Gamma'}$. Then for $\tau\to 0$, one has point-wisely that
        \begin{align}
       & \bar{u}^{(\tau)}=\bar{u}+o(1),~ \bar{v}^{(\tau)}=\bar{v}+o(1),~    \bar{w}^{(\tau)}=\bar{w}+o(1),\nonumber\\
       & \bar{E}^{(\tau)}=\bar{E}+o(1),~ \frac{\mathrm{d}\bar{\rho}^{(\tau)}}{\mathrm{d}\bar{\mu}}=\tau\frac{\mathrm{d}\bar{\rho}}{\mathrm{d}\bar{\mu}}+o(\tau), ~\bar{w}_p^{(\tau)}=\bar{w}_p+o(1), \label{eq418}
        \end{align}
where $\bar{U}=(\bar{\varrho},\bar{u},\bar{v},\bar{w},\bar{E})^{\top}$ is the axisymmetric Radon measure solution to $Problem ~B^{\prime}$, and $\bar{w}_p^{(\tau)}=w_p/(\gamma p_\infty K^2)$ (see \eqref{eq436}).  
\end{theorem}

The meaning of convergence claimed in \eqref{eq418} will be specified in Section \ref{sec 4.5}.  

%Notably, under the assumption of axisymmetry for the cone, it is natural to address this problem in the cylindrical coordinates. When we consider the axisymmetric solutions, the problem reduces to a two-dimensional case, albeit with the consequence that the flux depends on the radius $r$, rendering the governing equations non-homogeneous.

\subsection{Formulation of problems in cylindrical coordinates}\label{sec 4.2}

Under the assumption of axisymmetry of the cone, it is natural to reformulate Problem A$^{\prime}$ and Problem B$^{\prime}$ in cylindrical coordinates $\left(x, r \doteq \sqrt{y^{2}+z^{2}}, \theta\right) $ of $\mathbb{R}^{3}$, with
$$x=x, \quad y=r\cos \theta, \quad z=r \sin \theta,\quad r>0~\text{and}~\theta \in[-\pi, \pi).$$ 

Let us first restate Problem A$^{\prime}$ in cylindrical coordinates. We denote $u_x=u$ as the axial velocity, $u_r=v\cos\theta+w\sin\theta$ as the radial velocity, and $u_{\theta}=w\cos\theta-v\sin\theta$ as the swirl velocity. Then the Euler equations (\ref{Euler_eq 3D})--(\ref{state_eq 3D}) can be rewritten as 
\begin{equation}\label{eq in cc}
    \left\{\begin{aligned}
&\partial_x(r \rho u_x)+\partial_r(r \rho u_r)+\partial_{\theta}(\rho u_{\theta})=0,\\
&\partial_x{\left[r\left(\rho u_x^{2}+p\right)\right]+\partial_r(r \rho u_x u_r)+\partial_{\theta}(r \rho u_x u_{\theta})=0,} \\
&\partial_x(r \rho u_x u_r)+\partial_r\left[r\left(\rho u_r^{2}+p\right)\right]+\partial_{\theta}(\rho u_r u_{\theta})-\left(\rho u_{\theta}^{2}+p\right)=0, \\
&\partial_x(r \rho u_x u_{\theta})+\partial_r(r \rho u_r u_{\theta})+\partial_{\theta}\left(\rho u_{\theta}^{2}+p\right)+\rho u_r u_{\theta}=0, \\
&\partial_x(r \rho u_x E)+\partial_r(r \rho u_r E)+\partial_{\theta}(r \rho u_{\theta} E)=0 ,
\end{aligned}\right.
\end{equation}
and
\begin{equation}\label{p in cc}
p=\frac{\gamma -1}{\gamma }\rho\cdot\big(E-\frac{1}{2}(u_{x}^2+u_{r}^2+u_{\theta}^2)\big).
\end{equation}

For simplicity of notations but without misunderstanding, we abbreviate hereafter $(u_x,u_r,u_{\theta})$ to $(u,v,w)$.

Owing to the feature of the cone, we can construct axisymmetric solutions to Problem A$^{\prime}$, namely those solutions  independent of $\theta$, with  $w\equiv0$. Then (\ref{eq in cc})--(\ref{p in cc}) can be reduced to
\begin{equation}\label{Euler eq axisymmetric 3D}
    \left\{\begin{aligned}
&\partial_x(r \rho u)+\partial_r(r \rho v)=0,\\
&\partial_x{\left[r\left(\rho u^{2}+p\right)\right]+\partial_r(r \rho u v)=0,} \\
&\partial_x(r \rho u v)+\partial_r\left[r\left(\rho v^{2}+p\right)\right]-p=0, \\
&\partial_x(r \rho u E)+\partial_r(r \rho v E)=0,
\end{aligned}\right.
\end{equation}
and
\begin{equation}\label{state eq axisymmetric 3D}
    p=\frac{(\gamma-1)\rho}{\gamma}[E-\frac{1}{2}(u^2+v^2)].
\end{equation}
The corresponding fluid domain (\ref{domaineq}) and its boundary of (\ref{boundaryeq}) can be expressed as
\begin{equation*}
    \Omega_{as} \doteq \left\{(x, r):~x\ge0, ~ r>\tau f(x)\right\},\quad \Gamma_{as} \doteq\{(x, r):~ x\geq 0,~ r=\tau f(x)\}.
\end{equation*}
The boundary condition (\ref{slipcondition 3D}) on $\Gamma_{as}$ reduces to
\begin{equation}\label{slip boundary axisymmetric 3D}
    v=\tau u f^{\prime}(x),
\end{equation}
and the oncoming flow (\ref{oncoming flow 3D})--(\ref{oncoming E,p}) becomes
\begin{equation}\label{oncoming flow axisymmetric 3D}
    U_{\infty}=\left(\rho_{\infty}, u_{\infty}, 0, E_{\infty}\right)^{\top},\quad E_{\infty}=\frac{1}{2}u_{\infty}^2\Big(1+\frac{2\tau^2}{(\gamma -1)K^2}\Big),\quad 
        p_{\infty}=\frac{\rho_{\infty}u_{\infty}^2\tau^2}{\gamma K^2}.
\end{equation}

Similarly, for Problem B$^{\prime}$, Eqs. (\ref{hsd govern eq 3D})--(\ref{hsd state eq 3D}) in the  scaled fluid region $$\Omega'_{as} \doteq \{(\bar{x}, \bar{r}):~\bar{x}\ge0,~\bar{r}> f(\bar{x})\}$$ are simplified to
\begin{equation}\label{hsd govern axisymmetric 3D}
    \left\{\begin{aligned}
    &\partial_{\bar{x}}(\bar{r} \bar{\rho} )+\partial_{\bar{r}}(\bar{r} \bar{\rho} \bar{v})=0,\\
&\partial_{\bar{x}}{\left[\bar{r}\left(\bar{\rho} \bar{u}+\bar{p}\right)\right]+\partial_{\bar{r}}(\bar{r}\bar{\rho} u \bar{v})=0,} \\
&\partial_{\bar{x}}(\bar{r} \bar{\rho}  \bar{v})+\partial_{\bar{r}}\left[\bar{r}\left(\bar{\rho} \bar{v}^{2}+\bar{p}\right)\right]-\bar{p}=0, \\
&\partial_{\bar{x}}(\bar{r} \bar{\rho} E)+\partial_{\bar{r}}(\bar{r} \bar{\rho} \bar{v} E)=0,
\end{aligned}\right.
\end{equation}
and
\begin{equation}\label{hsd state axisymmetric 3D}
    \bar{p}=\frac{\gamma-1}{2\gamma}(\bar{E}-2\bar{u}-\bar{v}^2).
\end{equation}
The boundary condition (\ref{slip boundary_hsd 3D}) on the scaled boundary 
 $${\Gamma}'_{as} \doteq\{(\bar{x}, \bar{r}): \bar{x}\geq 0, ~ \bar{r}= f(\bar{x})\} $$ converts to
\begin{equation}\label{hsd bounadary axisymmetric 3D}
\bar{v}=f^{\prime}(\bar{x}),
\end{equation}
and the oncoming flow (\ref{oncoming flow_hsd 3D}) becomes
\begin{equation}\label{hsd oncoming axisymmetric 3D}
        \bar{U}_{\infty}=(1, 0, 0, \bar{E}_{\infty})^{\top},\quad \bar{E}_\infty=\frac{2}{(\gamma-1) K^2},\quad \bar{p}_\infty=\frac{1}{\gamma K^2}.
\end{equation}

\subsection{An axisymmetric Radon measure solution to Problem A$^{\prime}$}\label{sec 4.3}
We use ${\mathcal{L}}_{as}^{2}$ to denote the axisymmetric case of the Lebesgue measure $\mathcal{L}^{3}$ (the two differ by a factor $2\pi$), which is different from the standard Lebesgue measure $\mathcal{L}^{2}$. The Radon measure ${\mathcal{L}}^{2}_{as}$ on $\mathbb{R}^2_+\doteq\{(x,r): x\in \mathbb{R},~r\ge0\}$ is defined by
\begin{equation}\label{circular Radon measure}
\left\langle {\mathcal{L}}_{as}^{2}, \phi\right\rangle \doteq \int_{\Omega_{as}} \phi(x, r) r \mathrm{d} x\mathrm{d} r, \quad \forall \phi \in C_{c}(\mathbb{R}^2_+).
\end{equation}
Suppose that $L$ is a $C^1$ curve on $\mathbb{R}^2_+$, represented in parametric form as $(x(t), r(t))$ for $t \in[0, T)$, and $w_{L}(t)$ is a continuous function. Then the weighted circular Dirac measure can be defined as:
\begin{equation}\label{circular Dirac measure}
\left\langle w_{L} \delta_{L_{as}}, \phi\right\rangle \doteq \int_{0}^{T} r(t) w_{L}(t) \phi(x(t), r(t)) \sqrt{x^{\prime}(t)^{2}+r^{\prime}(t)^{2}}~ \mathrm{d} t, \quad \forall \phi \in C_{c}(\mathbb{R}^2_+).
\end{equation}
It is obvious that we have the following Radon-Nikodym derivatives: 
\begin{equation*}
\frac{\mathrm{d} \mathcal{L}_{as}^{2}}{\mathrm{d} \mathcal{L}^{2}}=r, \quad \frac{\mathrm{d} \delta_{L_{as}}}{\mathrm{d} \delta_{L}}=r.
\end{equation*}

We now raise the definition of an axisymmetric Radon measure solution to the Problem A$^{\prime}$ in cylindrical coordinates.
\begin{definition}
    For fixed $\gamma> 1$, let $\varrho,~\wp,~{\wp}'$ be nonnegative Radon measures, and $m_{a}^{j},~ m_{b}^{j}\, (j=0,1,2,3)$ be signed Radon measures on $\overline{\Omega}_{as}$, and $w_p(x)$ be a nonnegative function belonging to $L_{loc}^1([0,\infty))$. We call $(\varrho ,u,v,E  )$ a Radon
measure solution to $(\ref{Euler eq axisymmetric 3D})$--$(\ref{oncoming flow axisymmetric 3D})$, provided that

    $(i) [$linear relaxation$]$ for any $\phi\in C_c^1(\mathbb{R}^2_{+})$, there are
   \begin{equation}\label{3.25}
   \begin{aligned}
        \left\langle m_{a}^{0}, \partial_x\phi\right\rangle+\left\langle m_{b}^{0}, \partial_r\phi\right\rangle &+\int_{0}^{\infty} r \rho_{\infty} u_{\infty} \phi(0, r) \mathrm{d} r=0,\\
        \left\langle m_{a}^{1}+\wp, \partial_x\phi\right\rangle+\left\langle m_{b}^{1}, \partial_r\phi\right\rangle+\left\langle w_{p}(x) n_{1} \delta_{\Gamma_{as}},\phi\right\rangle&+\int_{0}^{\infty} r\left(\rho_{\infty} u_{\infty}^{2}+p_{\infty}\right) \phi(0, r) \mathrm{d}r=0,\\
       \left\langle m_{a}^{2}, \partial_x\phi \right\rangle+\left\langle m_{b}^{2}+\wp, \partial_r\phi \right\rangle&+\langle{\wp}', \phi\rangle+\left\langle w_{p}(x) n_{2} \delta_{\Gamma_{as}}, \phi\right\rangle =0,\\  
        \left\langle m_{a}^{3}, \partial_x\phi\right\rangle+\left\langle m_{b}^{3}, \partial_r\phi\right\rangle&+\int_{0}^{\infty} r \rho_{\infty} u_{\infty} E_{\infty} \phi(0, r) \mathrm{d} r=0,  
        \end{aligned}
    \end{equation}
    where $\mathbf{n}=(n_1,n_2)\doteq(-\tau f'(x),1)/\sqrt{1+\tau^2f'(x)^2}$ denotes the unit  normal vector on $\Gamma_{as}$ pointing to $\Omega_{as}$;
      
      $(ii) [$nonlinear constraints$]$
     there are   $\wp \ll \varrho, ~{\wp}' \ll \varrho,$ and $~(m_{a}^{j}, ~m_{b}^{j}) \ll$ $\varrho~(j=0,1,2,3)$, with the Radon-Nikodym derivatives 
   \begin{equation}\label{3.26}
       u=\frac{\mathrm{d}m_a^0}{\mathrm{d}\varrho}
       \quad \text{and}\quad   v=\frac{\mathrm{d}m_b^0}{\mathrm{d}\varrho}
   \end{equation}
   satisfying $\varrho$-a.e. that
   \begin{equation}
       \label{3.27}
       uu=\frac{\mathrm{d}m_a^1}{\mathrm{d}\varrho  } 
       ,\quad uv=\frac{\mathrm{d}m_b^1}{\mathrm{d}\varrho  }
            =\frac{\mathrm{d}m_a^2}{\mathrm{d}\varrho  },\quad
        vv=\frac{\mathrm{d}m_b^2}{\mathrm{d}\varrho},  
   \end{equation}   
    and there is a $\varrho$-a.e. function $E$ so that
      \begin{equation}
       \label{3.28}
        Eu=\frac{\mathrm{d}m_a^3}{\mathrm{d}\varrho  },\quad
        Ev=\frac{\mathrm{d}m_b^3}{\mathrm{d}\varrho  };
   \end{equation}

     $(iii)[$state equation$]$ if $\varrho \ll {\mathcal{L}}_{as}^{2}$ with Radon-Nikodym derivative $\rho(x, r)$, $\wp\ll{\mathcal{L}}_{as}^{2}$ with Radon-Nikodym derivative $p(x, r)$, and $\wp'\ll{\mathcal{L}}_{as}^{2}$ with Radon-Nikodym derivative $p'(x, r)$ in a neighborhood of $(x, r)\in\Omega_{as}$, then ${\mathcal{L}}_{as}^{2}$-a.e. there are
        \begin{equation*}
    \rho=\frac{\mathrm{d}\varrho}{\mathrm{d}{\mathcal{L}}_{as}^{2}},\quad p=\frac{\mathrm{d}\wp}{\mathrm{d}{\mathcal{L}}_{as}^{2}}=\frac{\gamma-1}{\gamma} \rho \left(E-\frac{u^{2}+v^{2}}{2}\right), \quad p'=\frac{\mathrm{d}\wp'}{\mathrm{d}{\mathcal{L}}_{as}^{2}}=p.
    \end{equation*}
    Furthermore, the classical entropy condition is valid for discontinuities of functions $\rho, u, v, E$ in this case.
\end{definition}

For the existence of Radon measure solutions to Problem A$^{\prime}$, we have the following result.
\begin{theorem}\label{th3.2}
       Assume that $f(x)\in C^2([0,\infty))$ and satisfies
    \begin{equation}\label{existence condition problem A'}
    \begin{aligned}
        f(0)=0,&\quad f'(x)> 0,\\
        f(x)(1+\tau^2f'(x)^2)^{\frac{3}{2}}+\gamma K^2 f''(x)M(x)&+\sqrt{1+\tau^2f'(x)^2}~f(x)f'(x)^2\ge 0
    \end{aligned}
    \end{equation}
    for $x\ge0$, 
    where
    \begin{equation*}
       M(x)\doteq \int_0^{x}\frac{f(t)f'(t)}{\sqrt{1+\tau^2f'(t)^2}}\mathrm{d}t.
    \end{equation*}
    Then Problem A$^\prime$ admits an axisymmetric Radon measure solution provided by $(\ref{solution u,v to Problem A'})$. In particular, we have for $x\ge0$ that 
    \begin{equation}\label{eq436}
       w_p(x)=\rho_{\infty}u_{\infty}^2\tau^2\Big(\frac{1}{\gamma K^2}+\frac{f''(x)M(x)+\sqrt{1+\tau^2f'(x)^2}f(x)f'(x)^2}{f(x)(1+\tau^2f'(x)^2)^{\frac{3}{2}}}\Big).
    \end{equation}
\end{theorem}

To prove Theorem \ref{th3.2}, we construct a Radon measure solution to (\ref{Euler eq axisymmetric 3D})--(\ref{oncoming flow axisymmetric 3D}) given by measures of the form 
\begin{equation}\label{A' solution of the form_1}
m_{a}^{0} \doteq \rho_{\infty} u_{\infty} {\mathcal{L}}_{as}^{2}\mres\Omega_{as}+w_{a}^{0}(x) \delta_{\Gamma_{as}}, \quad m_{b}^{0} \doteq w_{b}^{0}(x) \delta_{\Gamma_{as}},
\end{equation}
\begin{equation}\label{A' solution of the form_2}
m_{a}^{1} \doteq \rho_{\infty} u_{\infty}^2{\mathcal{L}}_{as}^{2}\mres\Omega_{as}+w_{a}^{1}(x) \delta_{\Gamma_{as}}, \quad m_{b}^{1} \doteq w_{b}^{1}(x) \delta_{\Gamma_{as}}, 
\end{equation}
\begin{equation}
 \wp \doteq p_{\infty}{\mathcal{L}}_{as}^{2}\mres\Omega_{as}, \quad \wp' \doteq p_{\infty}{\mathcal{L}}_{as}^{2}\mres\Omega_{as},
\end{equation}
\begin{equation}\label{A' solution of the form_3}
m_{a}^{2} \doteq w_{a}^{2}(x) \delta_{\Gamma_{as}}=w_{b}^{1}(x) \delta_{\Gamma_{as}}, \quad m_{b}^{2} \doteq w_{b}^{2}(x) \delta_{\Gamma_{as}},
\end{equation}
\begin{equation}\label{A' solution of the form_4}
m_{a}^{3} \doteq \rho_{\infty}u_{\infty}E_{\infty} {\mathcal{L}}_{as}^{2}\mres\Omega_{as}+w_{a}^{3}(x) \delta_{\Gamma_{as}}, \quad m_{b}^{3} \doteq w_{b}^{3}(x) \delta_{\Gamma_{as}}, 
\end{equation}
where  $w_{a}^{j}(x), w_{b}^{j}(x)~(j=0,1,2,3)$  are undetermined $C^1$ functions called weights.

We now derive the ordinary differential equations obeyed by $w_a^j(x), w_b^j(x), j=0,1,2,3$. Substituting (\ref{A' solution of the form_1}) into (\ref{3.25})$_1$, through the definition of Radon measures (\ref{circular Radon measure})--(\ref{circular Dirac measure}), we derive
\begin{equation}\label{derive of w_ab0}
\begin{aligned} 
&\int _{\Omega_{as} }r\rho_{\infty} u_{\infty}\partial _x\phi \mathrm{d}x\mathrm{d}r+\int _0^{\infty}\tau f(x)\sqrt{1+\tau^2f'(x)^2}w_a^0(x)\partial _x\phi(x,\tau f(x))\mathrm{d}x\\
&\quad + \int _0^{\infty}\tau f(x)\sqrt{1+\tau^2f'(x)^2}w_b^0(x)\partial _r\phi(x,\tau f(x))\mathrm{d}x+\int _0^{\infty}r\rho_{\infty} u_{\infty}\phi(0,r)\mathrm{d}r=0.  
\end{aligned}
\end{equation}
By Gauss-Green theorem, it holds that 
    \begin{equation}\label{gaussgreenf}
    \begin{aligned}
&\int _{\Omega_{as} }r\rho_{\infty} u_{\infty}\partial_x\phi(x,r)\,\mathrm{d}x\mathrm{d}r
=\int_{\partial \Omega_{as} }r\rho_{\infty} u_{\infty}\phi(x,r)\,\mathrm{d} r\\
=& \int_0^{\infty}\tau^2 f'(x) f(x)\rho_{\infty} u_{\infty}\phi (x,\tau f(x))\mathrm{d}x-\int_0^{\infty} r\rho_{\infty} u_{\infty}\phi(0,r)\mathrm{d}r.  
\end{aligned}
\end{equation}
Through the chain-rule $ \frac{\mathrm{d}\phi}{\mathrm{d}x}(x,\tau f(x))=\tau f'(x)\partial_r\phi(x,\tau f(x))+\partial_x\phi(x,\tau f(x)$, it holds 
 \begin{equation} \label{chainrule}
    \begin{aligned}
&\int _0^{\infty}\tau f(x)\sqrt{1+\tau^2f'(x)^2}w_b^0(x)\partial _r\phi(x,\tau f(x))\mathrm{d}x\\
=&\int _0^{\infty}\tau f(x)\sqrt{1+\tau^2f'(x)^2}w_b^0(x)\frac{\frac{\mathrm{d}\phi}{\mathrm{d}x}-\partial_x\phi  }{\tau f'(x)} \mathrm{d}x.
\end{aligned}
\end{equation}
From (\ref{gaussgreenf})--(\ref{chainrule}) and expression of the boundary $r=\tau f(x)$, for any test function $\phi\in C_{c}^1(\overline{\Omega}_{as})$, (\ref{derive of w_ab0}) can be rewritten as 
\begin{equation*}
\begin{aligned}
&\int_0^{\infty}\phi(x,\tau f(x))\Big[\tau  f(x)f'(x)\rho_{\infty}u_{\infty}-\frac{\mathrm{d}(\frac{f(x)\sqrt[]{1+\tau^2f'(x)^2}~w_b^0(x)}{f'(x)})}{\mathrm{d}x } \Big]\mathrm{d}x\\
&\quad +\int_0^{\infty}\tau f(x) \sqrt[]{1+\tau^2f'(x)^2}\,\partial_x\phi(x,\tau f(x))\Big[w_a^0(x)-\frac{w_b^0(x)}{\tau f'(x)}\Big]\, \mathrm{d}x=0. 
\end{aligned}
\end{equation*}
Therefore, for any $x>0$, the weight function $w_a^0(x)$ and $w_b^0(x)$ satisfy
\begin{equation}\label{ODEof w_ab0}
\left\{\begin{aligned}
         &\tau  f(x)f'(x)\rho_{\infty}u_{\infty}-\frac{\mathrm{d}\Big(\big( f(x)\sqrt[]{1+\tau^2f'(x)^2}~w_a^0(x)\big)/f'(x)\Big)}{\mathrm{d}x } =0,\\
&w_a^0(x)-\frac{w_b^0(x)}{\tau f'(x)}=0. 
\end{aligned}\right.
\end{equation}
By direct calculation, one gets
    \begin{equation}
        \label{3.35}
        \begin{aligned}
        w_{a}^{0}(x) &=\frac{\tau^2\rho_{\infty}u_{\infty} \int_0^{x}f(t)f'(t)\mathrm{d}t}{\tau f(x)\sqrt{1+\tau^2f'(x)^{2}}}=\frac{\tau\rho_{\infty}u_{\infty} f(x)}{2\sqrt{1+\tau^2f'(x)^{2}}},\\
        w_{b}^{0}(x)&=\tau f'(x) w_{a}^{0}(x)=\frac{\tau\rho_{\infty}u_{\infty}f(x)f'(x)}{2\sqrt{1+\tau^2f'(x)^{2}}}.
        \end{aligned}
    \end{equation}
%From the continuity of $w_{a}^{0}(x)$ and $w_{b}^{0}(x)$, combining with $f(0)=0$, we obtain
%$$w_{a}^{0}(0)=\lim _{x \rightarrow 0+} \frac{\tau f(x)}{\sqrt{1+\tau^2f^{\prime}(x)^{2}}}=0,\quad w_{b}^{0}(0)=\lim _{x \rightarrow 0+}\frac{\rho_{\infty}u_{\infty}\tau^2 f(x)}{2\sqrt{1+\tau^2f'(x)^{2}}}=0 .$$ 

Similarly, by substituting (\ref{A' solution of the form_2})--(\ref{A' solution of the form_4}) into (\ref{3.25}), we derive that 
\begin{equation}\label{axisymmetric ODE 3D}
     \left\{\begin{aligned}
         &\frac{\mathrm{d}\left(\tau f(x) \sqrt{1+\tau^2f'(x)^{2}}~w_{a}^{1}(x)\right)}{\mathrm{d} x}=\tau^2f(x)f'(x)(\rho_{\infty}u_{\infty}^2+p_{\infty}-w_p(x)), \\
         &\frac{\mathrm{d}\left(\tau f(x) \sqrt{1+\tau^2f'(x)^{2}}~w_{a}^{2}(x)\right)}{\mathrm{d} x}=\tau f(x)(w_{p}(x)-p_{\infty}), \\
         &\frac{\mathrm{d}\left(\tau f(x) \sqrt{1+\tau^2f'(x)^{2}}~w_{a}^{3}(x)\right)}{\mathrm{d} x}=\tau^2\rho_{\infty} u_{\infty} E_{\infty} f(x) f^{\prime}(x),
        \end{aligned}\right.
\end{equation}
and
\begin{equation}\label{axisymmetric ODE initial 3D}
w_{b}^{j}(x)=\tau f'(x) w_{a}^{j}(x), \quad j=1,2,3 . 
\end{equation}
Computing $(\ref{axisymmetric ODE 3D})_3$ and (\ref{axisymmetric ODE initial 3D}) directly, one yields that
\begin{equation}\label{3.36}
\begin{aligned}
w_{a}^{3}(x) &=\frac{\tau^2\rho_{\infty}u_{\infty}E_{\infty} \int_0^{x}f(t)f'(t)\mathrm{d}t}{\tau f(x)\sqrt{1+\tau^2f'(x)^{2}}}=\frac{\tau\rho_{\infty}u_{\infty}E_{\infty} f(x)}{2\sqrt{1+\tau^2f'(x)^{2}}}, \\
w_{b}^{3}(x)&=\tau f'(x) w_{a}^{3}(x)=\frac{\tau^2\rho_{\infty}u_{\infty}E_{\infty} f(x)f'(x)}{2\sqrt{1+\tau^2f'(x)^{2}}}.
\end{aligned}
\end{equation}

Noticing that we need to solve $w_p(x), w_a^1(x)$ and $w_a^2(x)$ from (\ref{axisymmetric ODE 3D})$_{1}$--(\ref{axisymmetric ODE 3D})$_{2}$. For $w_p(x)$, from (\ref{axisymmetric ODE 3D})$_{1}$, we obtain
\begin{equation}
w_p(x)=\rho_{\infty}u_{\infty}^2+p_{\infty}-\frac{1}{\tau^2f(x)f'(x)}\frac{\mathrm{d}(\tau f(x)\sqrt[]{1+\tau^2f'(x)^2}~w_a^1(x) ) }{\mathrm{d}x }.
\end{equation}
Now substituting it into (\ref{axisymmetric ODE 3D})$_2$, considering (\ref{axisymmetric ODE initial 3D}) and the relation $w_b^1(x)=w_a^2(x)$, one gets
\begin{equation}
\begin{aligned}
 & \tau f'(x)\frac{\mathrm{d}(\tau f(x)\sqrt[]{1+\tau^2f'(x)^2}~\tau f'(x)w_a^1(x) ) }{\mathrm{d}x}+\frac{\mathrm{d}(\tau f(x)\sqrt[]{1+\tau^2f'(x)^2}~w_a^1(x) ) }{\mathrm{d}x }\\
&\quad =\rho_{\infty}u_{\infty}^2\tau^2f(x)f'(x),  
\end{aligned}
\end{equation}
which can be simplified to 
\begin{equation}
\begin{aligned}
&(1+\tau^2 f'(x)^2)\frac{\mathrm{d}(\tau f(x)\sqrt[]{1+\tau^2f'(x)^2}~w_a^1(x) ) }{\mathrm{d}x}\\
&\qquad\quad+\tau^3 f(x)f'(x)f''(x)\sqrt{1+\tau^2f'(x)^2}~w_a^1(x)\\&\quad =\rho_{\infty}u_{\infty}^2\tau^2f(x)f'(x).
\end{aligned}
\end{equation}
 Then we divide the both sides by $\sqrt{1+\tau^2f'(x)^2}$ to obtain
 \begin{equation}\label{eq453}
\begin{aligned}
 & \sqrt{1+\tau^2f'(x)^2}\frac{\mathrm{d}(\tau f(x)\sqrt[]{1+\tau^2f'(x)^2}~w_a^1(x) ) }{\mathrm{d}x}+\tau^3 f(x)f'(x)f''(x)w_a^1(x)\\&\quad =\frac{\rho_{\infty}u_{\infty}^2\tau^2f(x)f'(x)}{\sqrt{1+\tau^2f'(x)^2}}. 
\end{aligned}
\end{equation}
If we set $A(x)=\tau f(x)\sqrt[]{1+\tau^2f'(x)^2}~w_a^1(x),~ B(x)=\sqrt{1+\tau^2f'(x)^2}$, then Eq. \eqref{eq453} reads 
\begin{equation}
    B(x)\frac{\mathrm{d}A(x)}{\mathrm{d}x}+A(x)\frac{\mathrm{d}B(x)}{\mathrm{d}x}=\frac{\rho_{\infty}u_{\infty}^2\tau^2f(x)f'(x)}{\sqrt{1+\tau^2f'(x)^2}}.
\end{equation}
Therefore, we have
\begin{equation}\label{4.52}
    \begin{aligned}  w_a^1(x)&=\frac{\tau\rho_{\infty}u_{\infty}^2M(x)}{ f(x)(1+\tau^2f'(x)^2)},\\
     w_b^1(x)&=\tau f'(x)w_a^1(x)=\frac{\tau^2\rho_{\infty}u_{\infty}^2f'(x)M(x)}{ f(x)(1+\tau^2f'(x)^2)},\\
     w_a^2(x)&= w_b^1(x)=\frac{\tau^2\rho_{\infty}u_{\infty}^2f'(x)M(x)}{ f(x)(1+\tau^2f'(x)^2)},\\
     w_b^2(x)&=\tau f'(x)w_a^2(x)=\frac{\tau^3\rho_{\infty}u_{\infty}^2f'(x)^2M(x)}{ f(x)(1+\tau^2f'(x)^2)},
    \end{aligned}
\end{equation}
and
\begin{equation}
    \label{4.53}   w_p(x)=p_{\infty}+\tau^2\rho_{\infty}u_{\infty}^2\frac{f''(x)M(x)+\sqrt{1+\tau^2f'(x)^2}f(x)f'(x)^2}{f(x)(1+\tau^2f'(x)^2)^{\frac{3}{2}}},
\end{equation}
where
\begin{equation}
    \label{3.40}
    M(x)\doteq \int_0^{x}\frac{f(t)f'(t)}{\sqrt{1+\tau^2f'(t)^2}}\mathrm{d}t.
\end{equation}
%From the continuity of $w_{a}^{1}(x)$,
%\begin{equation*}
%w_{a}^{1}(0)=\lim _{x \rightarrow 0+} w_{a}^{1}(x)=\lim _{x \rightarrow 0+} \frac{\tau\rho_{\infty}u_{\infty}^2M(x)}{ f(x)(1+\tau^2f'(x)^2)}=0 .
%\end{equation*}
%Similarly,
%\begin{equation*}
%    w_{a}^{2}(0)=\lim _{x \rightarrow 0+} w_{a}^{2}(x)=\lim _{x \rightarrow 0+} \frac{\tau^2\rho_{\infty}u_{\infty}^2f'(x)M(x)}{ f(x)(1+\tau^2f'(x)^2)}=0,
%\end{equation*}
%\begin{equation*}
%w_{p}(0)=\lim _{x \rightarrow 0+} w_{p}(x)=p_{\infty} .
%\end{equation*}

Noting (\ref{3.26})--(\ref{3.28}), from the related Radon-Nikodym derivatives we have 
\begin{equation*}
\left.u\right|_{\Gamma_{as}}=\frac{2 u_{\infty}M(x)}{f(x)^2 \sqrt{1+\tau^2f'(x)^{2}}},\left.\quad v\right|_{\Gamma_{as}}=\frac{2\tau u_{\infty} f^{\prime}(x) M(x)}{f(x)^2 \sqrt{1+\tau^2f^{\prime}(x)^{2}}},
\end{equation*}
and
\begin{equation*}
\left.E\right|_{\Gamma_{as}}=E_{\infty}.
\end{equation*}
Besides, the singular part of $\varrho$ is 
\begin{equation*}
    \frac{\tau \rho_{\infty}f(x)^3}{4 M(x)}\delta_{\Gamma_{as}}.
\end{equation*}
Finally, we show an axisymmetric Radon measure solution to Problem A$^{\prime}$ in cylindrical coordinates as following
\begin{equation}\label{solution u,v to Problem A'}
\begin{aligned}
u&=u_{\infty}\chi_{\Omega_{as}}+\frac{2 u_{\infty}M(x)}{f(x)^2 \sqrt{1+\tau^2f'(x)^{2}}} \chi_{\Gamma_{as}}, \quad v=\frac{2\tau u_{\infty} f^{\prime}(x) M(x)}{f(x)^2 \sqrt{1+\tau^2f^{\prime}(x)^{2}}} \chi_{\Gamma_{as}},\\
\varrho&=\rho_{\infty}{\mathcal{L}}_{as}^{2}\mres\Omega_{as}+\frac{\tau \rho_{\infty}f(x)^3}{4 M(x)} \delta_{\Gamma_{as}},
\quad 
E=E_{\infty} \chi_{\Omega_{as}}+E_{\infty} \chi_{\Gamma_{as}},
\end{aligned}
\end{equation}
with
$$
\begin{aligned}
& u(0+, 0)=\lim _{x \rightarrow 0+} \frac{2 u_{\infty}M(x)}{f(x)^2 \sqrt{1+\tau^2f'(x)^{2}}}=\frac{u_{\infty}}{1+\tau^2f^{\prime}\left(0_{+}\right)^{2}}, \\
& v(0+, 0)=\lim _{x \rightarrow 0+} \frac{2\tau u_{\infty} f^{\prime}(x) M(x)}{f(x)^2 \sqrt{1+\tau^2f^{\prime}(x)^{2}}}=\frac{\tau u_{\infty}f^{\prime}\left(0_{+}\right)}{1+\tau^2f^{\prime}\left(0_{+}\right)^{2}},\\
& \lim_{x\to 0+}\frac{\tau \rho_{\infty}f(x)^3}{4 M(x)}=0.\end{aligned}
$$
The proof of Theorem \ref{th3.2} is completed.

\subsection{An axisymmetric Radon measure solution to Problem B$^{\prime}$}\label{sec 4.4}

Let  $\bar{\mathcal{L}^{3}}$  be the Lebesgue measure in $(\bar{x},\bar{y},\bar{z})$-space. We use $\bar{\mathcal{L}}_{as}^{2}$ to denote the axisymmetric case of $\bar{\mathcal{L}^{3}}$,  and the two differ by a factor of $2\pi$. The Radon measure $\bar{\mathcal{L}}_{as}^{2}$ on the upper half-plane $\bar{\mathbb{R}}_{+}^{2} \doteq\{(\bar{x},\bar{r}): \bar{x} \in \mathbb{R},~ \bar{r} \geq 0\}$ is then defined as

\begin{equation}\label{circular Radon measure under scaling}
\left\langle \bar{\mathcal{L}}_{as}^{2}, \phi\right\rangle \doteq \int_{\bar{\mathbb{R}}_{+}^{2}} \phi(\bar{x}, \bar{r}) \bar{r} \mathrm{d} \bar{x}\mathrm{d} \bar{r}, \quad \forall \phi \in C_{c}\left(\bar{\mathbb{R}}_{+}^{2}\right).
\end{equation}

%\begin{remark}\rm 
%    Notably, as we mentioned in Section $\ref{sec 3.2}$, the Lebesgue measure $(\ref{circular Radon measure under scaling})$ is defined by neglecting scaling factor $\tau$. \hfill\qed
%\end{remark}

 Suppose that $L$ is a $C^1$ curve in $\overline{\Omega'_{as}}$. For $t \in[0, T)$, it can be parameterized as $(\bar{x}(t), \bar{r}(t))$. For a function  $\bar{w}_L(t)\in L^1_{\mathrm{loc}}([0, T])$, the Dirac measure supported on $L$ with the weight $\bar{w}_L$, which is singular with respect to $\bar{\mathcal{L}}_{as}^{2}$, can be defined by
    \begin{equation}\label{circular Dirac measure under scaling}
    \left\langle \bar{w}_{L} \bar{\delta}_{L_{as}}, \phi\right\rangle \doteq \int_{0}^{T} \bar{r}(t) \bar{w}_{L}(t) \phi(\bar{x}(t), \bar{r}(t)) \sqrt{\bar{x}^{\prime}(t)^{2}+\bar{r}^{\prime}(t)^{2}}\,\mathrm{d} t, \quad \forall\,\phi \in C_{c}\left(\bar{\mathbb{R}}_{+}^{2}\right).
    \end{equation}

\begin{definition}
       For $\gamma> 1$, let $\bar{\varrho},~\bar{\wp}, ~\bar{\wp}'$ be nonnegative Radon measures,~$\bar{m}_{a}^{j},~ \bar{m}_{b}^{j}~(j=0,1,2,3)$ be signed Radon measures on  $\overline{\Omega'_{as}}$, and $\bar{w}_p(\bar{x})\in L_{loc}^1([0,\infty) )$ be a nonnegative function. We call $(\bar{\varrho} ,\bar{u},\bar{v},\bar{E}  )$ a Radon measure solution to $(\ref{hsd govern axisymmetric 3D})$--$(\ref{hsd oncoming axisymmetric 3D})$, provided that
       
    $(i)[$linear relaxation$]$ for any $\phi\in C_c^1(\bar{\mathbb{R}}_+^2)$, there are   
            \begin{equation}\label{3.46}
        \begin{aligned}
        \left\langle \bar{m}_{a}^{0}, \partial_{\bar{x}}\phi\right\rangle+\left\langle \bar{m}_{b}^{0},\partial_{\bar{r}}\phi\right\rangle&+\int_{0}^{\infty}  \bar{r} \phi(0,\bar{r}) \mathrm{d} \bar{r}=0,\\
              \left\langle \bar{m}_{a}^{1}+\bar{\wp}, \partial_{\bar{x}}\phi\right\rangle+\left\langle \bar{m}_{b}^{1}, \partial_{\bar{r}}\phi\right\rangle&+\left\langle \bar{w}_{p}(\bar{x}) n'_{1} \bar{\delta}_{{\Gamma'_{as}}},\phi\right\rangle+\int_{0}^{\infty}  \bar{r}\bar{p}_{0} \phi(0, \bar{r}) \mathrm{d}\bar{r}=0,\\
       \left\langle \bar{m}_{a}^{2}, \partial_{\bar{x}}\phi \right\rangle+\left\langle \bar{m}_{b}^{2}+\bar{\wp}, \partial_{\bar{r}}\phi \right\rangle&+\langle\bar{\wp}', \phi\rangle+\left\langle \bar{w}_{p}(\bar{x}) n'_{2} \bar{\delta}_{{\Gamma'_{as}}}, \phi\right\rangle =0,  \\
       \left\langle \bar{m}_{a}^{3}, \partial_{\bar{x}}\phi\right\rangle+\left\langle \bar{m}_{b}^{3}, \partial_{\bar{r}}\phi\right\rangle&+\int_{0}^{\infty} \bar{r}\bar{E}_\infty \phi(0, \bar{r}) \mathrm{d} \bar{r}=0, 
    \end{aligned}
    \end{equation}
where $\mathbf{n}'=\left(n'_{1}, n'_{2}\right) \doteq \left(-f^{\prime}(\bar{x}), 1\right)\frac{1}{\sqrt{1+f^{\prime}(\bar{x})^{2}}}$ denotes the unit normal vector on ${\Gamma}'_{as}$ pointing into $\Omega'_{as}$, and $\bar{m}_{a}^{0}=\bar{\varrho}$;

$(ii)[$nonlinear constraints$]$  there are $\bar{\wp} \ll \bar{\varrho},~ \bar{\wp}' \ll \bar{\varrho},$ and $~\left(\bar{m}_{a}^{j}, \bar{m}_{b}^{j}\right) \ll \bar{\varrho}$$\,(j=0,1,2,3)$,  with corresponding Radon-Nikodym derivatives 
   \begin{equation}
       \label{3.47}     \bar{u}=\frac{\mathrm{d}\bar{m}_a^1}{\mathrm{d}\bar{\varrho}  }\quad \text{and} \quad   \bar{v}=\frac{\mathrm{d}\bar{m}_b^0}{\mathrm{d}\bar{\varrho}}=\frac{\mathrm{d}\bar{m}_a^2}{\mathrm{d}\bar{\varrho}}
   \end{equation}
   satisfying $\bar{\varrho}$-a.e. that
        \begin{equation} \label{3.48}  \bar{u}\bar{v}=\frac{\mathrm{d}\bar{m}_b^1}{\mathrm{d}\bar{\varrho}  }
            ,\quad      \bar{v}\bar{v}=\frac{\mathrm{d}\bar{m}_b^2}{\mathrm{d}\bar{\varrho}},
        \end{equation}
     and there is a $\bar{\varrho}$-a.e. function $\bar{E}$ so that
   \begin{equation}\label{3.49}    \bar{E}=\frac{\mathrm{d}\bar{m}_a^3}{\mathrm{d}\bar{\varrho}  },\quad    \bar{E}\bar{v}=\frac{\mathrm{d}\bar{m}_b^3}{\mathrm{d}\bar{\varrho}  };
   \end{equation}

   $(iii)[$state equation$]$ if $\bar{\varrho} \ll \bar{\mathcal{L}}_{as}^{2}$,$~\bar{\wp} \ll \bar{\mathcal{L}}_{as}^{2}$ and $\bar{\wp}' \ll \bar{\mathcal{L}}_{as}^{2}$ in a neighborhood of $(\bar{x},\bar{r})\in \Omega_{as}'$, then  $\bar{{\mathcal{L}}}_{as}^{2}$-a.e. it holds 
        \begin{equation*}
   \bar{\rho}=\frac{\mathrm{d}\bar{\varrho}}{\mathrm{d}\bar{\mathcal{L}}_{as}^{2}},\quad \bar{p}=\frac{\mathrm{d}\bar{\wp}}{\mathrm{d}\bar{\mathcal{L}}_{as}^{2}}=\frac{{\gamma}-1}{2{\gamma}} \bar{\rho} \left(\bar{E}-2\bar{u}-\bar{v}^2\right), \quad \bar{p}'=\frac{\mathrm{d}\bar{\wp}'}{\mathrm{d}\bar{\mathcal{L}}_{as}^{2}}= \bar{p}.
    \end{equation*}
    Furthermore, the classical entropy condition is valid for discontinuities of functions $\bar{U}=(\bar{\rho}, \bar{u},\bar{v},\bar{E})^{\top}$ in this case.  
\end{definition}

For the existence of a Radon measure solution to (\ref{hsd govern axisymmetric 3D})--(\ref{hsd oncoming axisymmetric 3D}), we have the following theorem.

\begin{theorem}\label{th3.3}
 Assume that $f(\bar{x})\in C^2([0,\infty))$, and for $\bar{x}\ge 0$, the following hold:
\begin{equation}\label{existence condition Problem B'}
\begin{aligned}
    f(0)=0,\quad f'(\bar{x})> 0,\\
    2f(\bar{x})+\gamma K^2(2f(\bar{x})f'(\bar{x})^2 + f(\bar{x})^2f''(\bar{x}))> 0.
\end{aligned}
\end{equation}
Then $(\ref{hsd govern axisymmetric 3D})$--$(\ref{hsd oncoming axisymmetric 3D})$ admits a Radon measure solution, which is corresponding to an axisymmetric Radon measure solution to Problem B$^{\prime}$. The associated $\bar{w}_p(\bar{x})$ is given by \eqref{hsd axis weight function w_p} below. 
\end{theorem}

We now prove Theorem \ref{th3.3} by constructing a Radon measure solution to (\ref{hsd govern axisymmetric 3D})--$(\ref{hsd oncoming axisymmetric 3D})$ given by measures of the form
        \begin{equation}\label{hsd solution of the form_1}
           \bar{\varrho}=\bar{m}_{a}^{0} \doteq  \bar{\mathcal{L}}_{as}^{2}\mres\Omega'_{as}+\bar{w}_{a}^{0}(\bar{x}) \bar{\delta}_{{\Gamma}'_{as}}, \quad \bar{m}_{b}^{0} \doteq \bar{w}_{b}^{0}(\bar{x}) \bar{\delta}_{{\Gamma}'_{as}},
      \end{equation}
      \begin{equation}\label{hsd solution of the form_2}
           \bar{m}_{a}^{1} \doteq \bar{w}_{a}^{1}(\bar{x}) \bar{\delta}_{{\Gamma}'_{as}}, \quad \bar{m}_{b}^{1} \doteq \bar{w}_{b}^{1}(\bar{x}) \bar{\delta}_{{\Gamma}'_{as}}, \quad \bar{\wp} \doteq \bar{p}_\infty \bar{\mathcal{L}}_{as}^{2}\mres\Omega'_{as} \quad \bar{\wp}' \doteq \bar{p}_\infty \bar{\mathcal{L}}_{as}^{2}\mres\Omega'_{as},
      \end{equation}
      \begin{equation}\label{hsd solution of the form_3}
           \bar{m}_{a}^{2} \doteq \bar{w}_{a}^{2}(\bar{x}) \bar{\delta}_{{\Gamma}'_{as}}=\bar{w}_{b}^{0}(\bar{x}) \bar{\delta}_{{\Gamma}'_{as}}, \quad \bar{m}_{b}^{2} \doteq \bar{w}_{b}^{2}(\bar{x}) \bar{\delta}_{{\Gamma}'_{as}},
      \end{equation}
      \begin{equation}\label{hsd solution of the form_4}
           \bar{m}_{a}^{3} \doteq \bar{E}_\infty \bar{\mathcal{L}}_{as}^{2}\mres\Omega'_{as}+\bar{w}_{a}^{3}(\bar{x}) \bar{\delta}_{{\Gamma}'_{as}}, \quad \bar{m}_{b}^{3} \doteq \bar{w}_{b}^{3}(\bar{x}) \bar{\delta}_{{\Gamma}'_{as}}, 
      \end{equation}
 where  $\bar{w}_{a}^{j}(\bar{x})$ and  $\bar{w}_{b}^{j}(\bar{x})\, (j=0,1,2,3)$  are undetermined $C^1$ functions.

 Substituting (\ref{hsd solution of the form_1})--(\ref{hsd solution of the form_4}) into (\ref{3.46}), we can derive the differential equations obeyed by these weight functions
    \begin{equation}
        \label{hsd axis ODE of weight functions 3D}
   \left\{\begin{aligned}
         &\frac{\mathrm{d}\left( f(\bar{x}) \sqrt{1+f'(\bar{x})^{2}}~\bar{w}_{a}^{0}(\bar{x})\right)}{\mathrm{d}{\bar{x}}}=f(\bar{x}) f^{\prime}(\bar{x}),\\
         &\frac{\mathrm{d}\left(f(\bar{x}) \sqrt{1+f'(\bar{x})^{2}}~\bar{w}_{a}^{1}(\bar{x})\right)}{\mathrm{d}{\bar{x}}}=f(\bar{x})f'(\bar{x})\big(\bar{p}_\infty-\bar{w}_p(\bar{x})\big), \\
        &\frac{\mathrm{d}\left( f(\bar{x}) \sqrt{1+f'(\bar{x})^{2}}~\bar{w}_{a}^{2}(\bar{x})\right)}{\mathrm{d}{\bar{x}}}=f(\bar{x})(\bar{w}_{p}(\bar{x})-\bar{p}_\infty), \\
         &\frac{\mathrm{d}\left( f(\bar{x}) \sqrt{1+f'(\bar{x})^{2}}~\bar{w}_{a}^{3}(\bar{x})\right)}{\mathrm{d} x}=\bar{E}_\infty f(\bar{x}) f^{\prime}(\bar{x}),
        \end{aligned}\right.
    \end{equation}
and
    \begin{equation}\label{hsd axis weight functions initial 3D}
\bar{w}_{b}^{j}(\bar{x})= f'(\bar{x}) \bar{w}_{a}^{j}(\bar{x}), \quad j=0,1,2,3 . 
\end{equation}

Similar calculations to those in Section \ref{sec 4.3} allow us to solve for $\bar{w}_{a}^{j}(\bar{x}), \bar{w}_{b}^{j}(\bar{x})\,(j=0,1,2,3)$ from Eqs. (\ref{hsd axis ODE of weight functions 3D})--(\ref{hsd axis weight functions initial 3D}), yielding the following results: 
\begin{equation}\label{hsd axis weight functions}
\begin{aligned} 
\bar{w}_{a}^{0}(\bar{x}) &=\frac{ \int_0^{\bar{x}}f(t)f'(t)\mathrm{d}t}{f(\bar{x})\sqrt{1+f'(\bar{x})^{2}}}=\frac{f(\bar{x})}{2\sqrt{1+f'(\bar{x})^{2}}}, \\
 \bar{w}_{b}^{0}(\bar{x})&=f'(\bar{x}) \bar{w}_{a}^{0}(\bar{x})=\frac{f(\bar{x})f'(\bar{x})}{2\sqrt{1+f'(\bar{x})^{2}}},\\
 \bar{w}_a^1(\bar{x})&=\frac{-f(\bar{x})^2f'(\bar{x})^2+\int_0^{\bar{x}}f(t)^2f'(t)f''(t)\mathrm{d}t}{ 2f(\bar{x})\sqrt{1+f'(\bar{x})^2}},\\
  \bar{w}_b^1(\bar{x})&=\frac{-f(\bar{x})^2f'(\bar{x})^3+f'(\bar{x})\int_0^{\bar{x}}f(t)^2f'(t)f''(t)\mathrm{d}t}{ 2f(\bar{x})\sqrt{1+f'(\bar{x})^2}},\\
  \end{aligned}
  \end{equation}
  as well as 
  \begin{equation}
\begin{aligned}
\bar{w}_a^2(\bar{x})&=\bar{w}_b^0(\bar{x})=\frac{f(\bar{x})f'(\bar{x})}{2\sqrt{1+f'(\bar{x})^{2}}},\\
\bar{w}_b^2(\bar{x})&=f'(\bar{x}) \bar{w}_a^2(\bar{x})=\frac{f(\bar{x})f'(\bar{x})^2}{2\sqrt{1+f'(\bar{x})^{2}}},\\
\bar{w}_{a}^{3}(\bar{x})&=\frac{\bar{E}_\infty \int_0^{\bar{x}}f(t)f'(t)\mathrm{d}t}{ f(\bar{x})\sqrt{1+f'(\bar{x})^{2}}}=\frac{\bar{E}_\infty f(\bar{x})}{2\sqrt{1+f'(\bar{x})^{2}}},\\
 \bar{w}_{b}^{3}(\bar{x})&=f'(\bar{x}) \bar{w}_{a}^{3}(\bar{x})=\frac{\bar{E}_\infty f(\bar{x})f'(\bar{x})}{2\sqrt{1+f'(\bar{x})^{2}}},
\end{aligned}
\end{equation}
and
\begin{equation}\label{hsd axis weight function w_p}
\bar{w}_p(\bar{x})=\bar{p}_\infty+\frac{1}{f(\bar{x})}\frac{\mathrm{d}\big(\frac{f(\bar{x})^2f'(\bar{x})}{2} \big)}{\mathrm{d}x }=\bar{p}_\infty+f'(\bar{x})^2+\frac{f(\bar{x})f''(\bar{x})}{2}.
\end{equation}
Then from  the constraints of  Radon-Nikodym derivatives, we have a Radon measure solution
\begin{equation}\label{a axis Radon measure solution  to B'}
\begin{aligned}
\bar{u}&=\Big(-f'(\bar{x})^2+\frac{ \int_0^{\bar{x}}f(t)^2f'(t)f''(t)\mathrm{d}t}{f(\bar{x})^2}\Big) \chi_{\Gamma_{as}}, \quad \bar{v}=f'(\bar{x}) \chi_{\Gamma_{as}},\\
\bar{\varrho}&=1\bar{\mathcal{L}}_{as}^{2}\mres\Omega'_{as}+\frac{f(\bar{x})}{2 \sqrt{1+f'(\bar{x})^2}} \bar{\delta}_{\Gamma'_{as}},\quad \bar{E}=\bar{E}_{\infty} \chi_{\Omega'_{as}}+\bar{E}_{\infty} \chi_{\Gamma'_{as}} 
\end{aligned}
\end{equation}
to (\ref{hsd govern axisymmetric 3D})--(\ref{hsd oncoming axisymmetric 3D}).  It is easy to extend these expressions to be valid for  $x=0$. The proof of Theorem \ref{th3.3} is completed.

\subsection{Comparison between Radon measure solutions for the two problems}\label{sec 4.5}
%We first clarify that, the $u$ and $\bar{u}$ that appear in (\ref{solution u,v to Problem A'}) and (\ref{a axis Radon measure solution u,v to B'}) are the horizontal velocity components $u_x$ and $\bar{u}_{\bar{x}}$ in the Cartesian coordinates, and the $v$ and $\bar{v}$ in the cylindrical coordinates represent the tangential velocity components $u_r$ and $\bar{u}_{\bar{r}}$ in the Cartesian coordinates.

%For Problem A$^{\prime}$, according to the cylindrical coordinate transformation, we can obtain
%\begin{equation*}
%   u=u_x,\quad v=u_r\cos{\theta}, \quad w=u_r\sin{\theta},
%\end{equation*}
%where
%$$\theta=\arctan {\frac{z}{y}}.$$
%Similarly, for Problem B$^{\prime}$,
%\begin{equation*}
%\bar{u}=\bar{u}_{\bar{x}},\quad \bar{v}=\bar{u}_{\bar{r}}\cos{\bar{\theta}},
%\quad \bar{w}=\bar{u}_{\bar{r}}\sin{\bar{\theta}},
%\end{equation*}
%where
%$$\bar{\theta}=\arctan{\frac{\bar{z}}{\bar{y}}}.$$
%From non-dimensional transformation $y=\tau \bar{y},~ z=\tau\bar{z}$, we can obtain
%$$\cos{\bar{\theta}}=\cos{\theta},\quad \sin{\bar{\theta}}=\sin{\theta}.$$
Technically, to establish  Theorem \ref{main theorem 3D}, the comparison between axisymmetric Radon measure solutions to Problem A$^{\prime}$ and Problem B$^{\prime}$ can be mapped to the cylindrical coordinates $(\bar{x}, \bar{r})$, which indicates that a comparison between $(\bar{\rho}^{(\tau)},\bar{u}^{(\tau)},\bar{v}^{(\tau)},\bar{E}^{(\tau)})$, the Radon measure solution $(\varrho,u,v, E)$ (see (\ref{solution u,v to Problem A'})) under the scaling \eqref{scaling 3D}, and $(\bar{\varrho},\bar{u},\bar{v},\bar{E})$ (see (\ref{a axis Radon measure solution  to B'})) suffices. 

Following Section \ref{sec2.3}, the proof is divided into two steps.

\textit{Step} 1. {\em Compare the functions $\bar{u}^{(\tau)},\bar{v}^{(\tau)},\bar{E}^{(\tau)}$ and $\bar{u},\bar{v},\bar{E}$.} It is easy to justify that the solutions away from the boundary $\Gamma_{as}'$ are constants independent of $\tau$. Here we establish the convergence of the solutions supported on the boundary $\Gamma_{as}'$ when slenderness $\tau\to 0$. 

From (\ref{solution u,v to Problem A'}) and (\ref{scaling 3D}), it is seen that 
\begin{equation*}
    \bar{u}^{(\tau)}|_{\Gamma_{as}'}=\frac{u|_{\Gamma_{as}}-u_{\infty}}{\tau^2u_{\infty}}=\frac{2{\int_{0}^{\bar{x}}\frac{ f(t)f'(t)}{\sqrt{1+\tau^2f'(t)^2}}\mathrm{d}t}-f(\bar{x})^2\sqrt{1+\tau^2f'(\bar{x})^2}}{\tau^2f(\bar{x})^2\sqrt{1+\tau^2 f'(\bar{x})^2}}.
\end{equation*}
Taking $\tau \to 0$, we have
\begin{equation}\label{3.65}
\begin{aligned}
    \lim_{\tau\to 0} \bar{u}^{(\tau)}|_{\Gamma_{as}'}=& \lim_{\tau \to 0} \frac{2{\int_{0}^{\bar{x}}\frac{ f(t)f'(t)}{\sqrt{1+\tau^2f'(t)^2}}\mathrm{d}t}-f(\bar{x})^2\sqrt{1+\tau^2f'(\bar{x})^2}}{\tau^2f(\bar{x})^2\sqrt{1+\tau^2 f'(\bar{x})^2}}\\  
    =&\lim_{\tau \to 0} \frac{2{\int_{0}^{\bar{x}}\frac{ f(t)f'(t)}{\sqrt{1+\tau^2f'(t)^2}}\mathrm{d}t}-f(\bar{x})^2\sqrt{1+\tau^2f'(\bar{x})^2}}{\tau^2f(\bar{x})^2}.
      \end{aligned}
\end{equation}
By L'Hospital rule, one infers that
\begin{align*}
         &\displaystyle  \lim_{\tau \to 0} \frac{2{\int_{0}^{\bar{x}}\frac{-\tau f(t)f'(t)^3}{(1+\tau^2f'(t)^2)^{\frac{3}{2}}}\mathrm{d}t}-\frac{\tau f(\bar{x})^2f'(\bar{x})^2}{\sqrt{1+\tau^2f'(\bar{x})^2}}}{2\tau f(\bar{x})^2}\\ 
         =&\lim_{\tau \to 0} \frac{{\int_{0}^{\bar{x}} - 2f(t)f'(t)^3\mathrm{d}t}- f(\bar{x})^2f'(\bar{x})^2}{2 f(\bar{x})^2}.
 \end{align*}
Moreover, thanks to the integration-by-part, we get
     \begin{equation*}
     \begin{aligned}
\int_{0}^{\bar{x}} - 2f(t)f'(t)^3\mathrm{d}t=-f(\bar{x})^2f'(\bar{x})^2+2\int_{0}^{\bar{x}}f(t)^2f'(t)f''(t)\mathrm{d}t.
\end{aligned}
\end{equation*}
So (\ref{3.65}) can be rewritten as
     \begin{equation}   \label{3.66}
        \lim_{\tau\to 0} \bar{u}^{(\tau)}|_{\Gamma_{as}'}
=-f'(\bar{x})^2+\frac{{\int_{0}^{\bar{x}}f(t)^2f'(t)f''(t)\mathrm{d}t}}{f(\bar{x})^2}=\bar{u}|_{\Gamma_{as}'}.
\end{equation}

Similarly, for velocity along $\bar{r}$-axis, one has
\begin{equation*}
    \bar{v}^{(\tau)}|_{\Gamma'}=\frac{v|_{\Gamma_{as}}}{\tau u_{\infty}},
\end{equation*}
 and as $\tau \to 0$,
\begin{equation}\label{3.67}
\begin{aligned}
        \lim_{\tau\to 0}\bar{v}^{(\tau)}|_{\Gamma_{as}'}
        =&\lim_{\tau \to 0} \frac{2 f^{\prime}(\bar{x}) \int_{0}^{\bar{x}}\frac{ f(t)f'(t)}{\sqrt{1+\tau^2f'(t)^2}}\mathrm{d}t}{f(\bar{x})^2 \sqrt{1+\tau^2f^{\prime}(\bar{x})^{2}}}\\
        =&\frac{2 f^{\prime}(\bar{x}) \int_{0}^{\bar{x}}f(t)f'(t)\mathrm{d}t}{f(\bar{x})^2 }\\
        =& f'(\bar{x})
        =\bar{v}|_{\Gamma_{as}'}.
 \end{aligned}
 \end{equation}
Besides, from (\ref{scaling_E 3D}) and (\ref{solution u,v to Problem A'}), we acquire
\begin{align}
    \label{3.69}
   \lim_{\tau\to 0}\bar{E}^{(\tau)}|_{\Gamma_{as}'}=\lim_{\tau \to 0}\frac{2E|_{\Gamma}-u_{\infty}^2}{u_{\infty}^2\tau^2}=\bar{E}_\infty= \bar{E}|_{\Gamma'_{as}}.
\end{align}
The verification of $(\ref{eq418})_3$ is also straightforward. 

\textit{Step} 2. We compare the circular measures $\bar{\varrho}^{(\tau)}$ and $\bar{\varrho}$, following a similar procedure 	as in Section \ref{sec2.3}.

For any test function $\phi \in C_c(\overline{\Omega}_{as})$, we recall 
\begin{equation*}
\begin{aligned}
    \langle\varrho,~\phi\rangle&=\langle\rho_{\infty}
    {\mathcal{L}}_{as}^2\mres{\Omega_{as}}+w_{\rho}(x)\delta_{\Gamma_{as}},~\phi\rangle\\
 %   &=\langle\rho_{\infty}{\mathcal{L}}_{as}^2\mres{\Omega_{as}},~\phi\rangle+\langle w_{\rho}(x)\delta_{\Gamma_{as}},~\phi\rangle\\
    &=\int_{\Omega_{as}}\rho_{\infty}r\phi~\mathrm{d}x\mathrm{d}r +\int_0^{\infty}w_{\rho}(x)\tau f(x)\phi(x,\tau f(x))\sqrt{1+\tau^2f'(x)^2}~ \mathrm{d}x.
\end{aligned}
\end{equation*}
Under the scaling (\ref{scaling 3D}), it reads 
\begin{equation}\label{3.70}
 \langle\bar{\varrho}^{(\tau)},~\bar{\phi}^{(\tau)}\rangle=\tau\Big(\int_{\Omega'_{as}}\bar{r}\bar{\phi}^{(\tau)}~\mathrm{d}\bar{x}\mathrm{d}\bar{r}+\int_0^{\infty}\frac{w_{\rho}(\bar{x})}{\tau \rho_{\infty}}f(\bar{x})\bar{\phi}^{(\tau)}(\bar{x},f(\bar{x}))\sqrt{1+\tau^2f'(\bar{x})^2}\,\mathrm{d}\bar{x}\Big).
\end{equation}
Then, considering $\bar{\varrho}$ given by \eqref{a axis Radon measure solution  to B'}, thanks to (\ref{circular Radon measure under scaling})--(\ref{circular Dirac measure under scaling}), for any test function $\bar{\phi}\in C_c(\overline{\Omega'}_{as})$, we have 
\begin{equation}\label{3.71}
\langle\bar{\varrho},~\bar{\phi}\rangle=\int_{\Omega'_{as}}\bar{r}\bar{\phi}~\mathrm{d}\bar{x}\mathrm{d}\bar{r} +\int_0^{\infty}\bar{w}_{\rho}(\bar{x})f(\bar{x})\bar{\phi}(\bar{x}, f(\bar{x}))\sqrt{1+f'(\bar{x})^2}~ \mathrm{d}\bar{x},
\end{equation}
where $\bar{w}_{\rho}(\bar{x})$ is $f(\bar{x})/(2 \sqrt{1+f'(\bar{x})^2})$, the singular part of $\bar{\varrho}$ (see \eqref{a axis Radon measure solution  to B'}). 

By ignoring $\tau$ in (\ref{3.70}), the mass of gas $\bar{\varrho}^{(\tau)}$ per unit $\mathrm{d}\bar{x}\mathrm{d}\bar{r}$  on $\Omega'_{as}$ is $\bar{r}$, which coincides with the mass of gas $\bar{\varrho}$ per unit $\mathrm{d}\bar{x}\mathrm{d}\bar{r}$ on $\Omega'_{as}$. The mass of gas $\bar{\varrho}^{(\tau)}$ per unit $\mathrm{d}\bar{x}$ on $\Gamma'_{as}$ is
\begin{equation*}
    \frac{w_{\rho}(\bar{x})}{\tau \rho_{\infty}}f(\bar{x})\sqrt{1+\tau^2f'(\bar{x})^2}=\frac{f(\bar{x})^4}{4\int_0^{\bar{x}}\frac{f(t)f'(t)}{\sqrt{1+\tau^2f'(t)^2}}\mathrm{d}t}\sqrt{1+\tau^2f'(\bar{x})^2},
\end{equation*}
compared with the mass of gas $\bar{\varrho}$ per unit $\mathrm{d}\bar{x}$ on $\Gamma'_{as}$, i.e., 
\begin{equation*}
    \bar{w}_{\rho}(\bar{x})f(\bar{x})\sqrt{1+f'(\bar{x})^2}=\frac{f(\bar{x})^2}{2},
\end{equation*}
then as $\tau \to 0$,  one gets
\begin{equation}
    \lim_{\tau\to 0}\frac{\frac{w_{\rho}(\bar{x})}{\tau \rho_{\infty}}\sqrt{1+\tau^2f'(\bar{x})^2}}{ \bar{w}_{\rho}(\bar{x})\sqrt{1+f'(\bar{x})^2}}=1.
\end{equation} 
This is what $(\ref{eq418})_2$ means.  The proof of Theorem \ref{thm41} is completed.

\appendix
\section{Mathematical structure of system (\ref{hyper sd eq})--(\ref{state eq for hsd})}\label{appendix a}
In this appendix, we are concerned with the mathematical structure of the new hypersonic small-disturbance equations (\ref{hyper sd eq})--(\ref{state eq for hsd}), which can be regraded as a model that is difficulty than the unsteady Euler system for one-dimensional non-isentropic compressible flow, but easier than the steady Euler system for two-dimensional non-isentropic supersonic flow.    

Substituting the equation of state (\ref{state eq for hsd}) into (\ref{hyper sd eq}), we obtain the following new hypersonic small-disturbance system:
    \begin{equation}\label{0.1}
        \left\{\begin{array}{l} 
    \partial _{\bar{x}}\bar{\rho}+\partial _{\bar{y}}(\bar{\rho} \bar{v})=0,\\  
    \partial _{\bar{x}}(\bar{\rho}\bar{u}+\frac{\gamma-1}{2\gamma}\bar{\rho}(\bar{E}-2\bar{u}-\bar{v}^2))+\partial _{\bar{y}}(\bar{\rho} \bar{u}\bar{v})=0,\\
     \partial _{\bar{x}}(\bar{\rho}\bar{v})+\partial _{\bar{y}}(\bar{\rho}\bar{v}^2+\frac{\gamma-1}{2\gamma}\bar{\rho}(\bar{E}-2\bar{u}-\bar{v}^2))=0,\\
             \partial _{\bar{x}}(\bar{\rho}\bar{E})+\partial _{\bar{y}}(\bar{\rho}\bar{v}\bar{E})=0,
        \end{array}\right. 
    \end{equation}
which can be written in the general form of conservation laws:
\begin{equation}\label{0.2}
W(\bar{U})_{\bar{x}}+H(\bar{U})_{\bar{y}}=0,
\end{equation}
where $\bar{U}=(\bar{\rho}, \bar{u}, \bar{v}, \bar{E})^\top$, and 
\begin{equation*}
   W(\bar{U})=(\bar{\rho},~\bar{\rho}\bar{u}+\frac{\gamma-1}{2\gamma}\bar{\rho}(\bar{E}-2\bar{u}-\bar{v}^2),~\bar{\rho}\bar{v},~\bar{\rho}\bar{E})^\top,
\end{equation*}
\begin{equation*}
    H(\bar{U})=(\bar{\rho}\bar{v},~\bar{\rho}\bar{u}\bar{v},~\bar{\rho}\bar{v}^2+\frac{\gamma-1}{2\gamma}\bar{\rho}(\bar{E}-2\bar{u}-\bar{v}^2),~\bar{\rho}\bar{v}\bar{E})^\top. 
\end{equation*}

For a smooth solution $\bar{U}(\bar{x},\bar{y})$, Eq. (\ref{0.2}) is equivalent to
\begin{equation*}   \nabla_{\bar{U}}W(\bar{U})\bar{U}_{\bar{x}}+\nabla_{\bar{U}}H(\bar{U})\bar{U}_{\bar{y}}=0.
\end{equation*}
The eigenvalues of (\ref{0.1}) are the roots of the fourth-order polynomial of $\lambda$, provided that $\bar{\rho}>0$:
\begin{equation}\label{0.3}
    \det(\lambda \nabla_{\bar{U}}W(\bar{U})-\nabla_{\bar{U}}H(\bar{U}))=0,
\end{equation}
where 
\begin{equation*}
\nabla_{\bar{U}}W(\bar{U})=\begin{pmatrix}
 1 & 0 & 0 & 0 \\
 \bar{u}+\frac{\gamma-1}{2\gamma}(\bar{E}-2\bar{u}-\bar{v}^2) & \frac{1}{\gamma }\bar{\rho} & \frac{1-\gamma}{\gamma}\bar{\rho}\bar{v} &\frac{\gamma -1}{2\gamma}\bar{\rho}   \\
 \bar{v} & 0  & \bar{\rho}   & 0 \\
  \bar{E} & 0  & 0 &\bar{\rho} 
\end{pmatrix},
\end{equation*}
and 
\begin{equation*}
    \nabla_{\bar{U}}H(\bar{U})=
    \begin{pmatrix}
 \bar{v}  & 0 & \bar{\rho}  & 0 \\
 \bar{u}\bar{v} & \bar{\rho}\bar{v}  & \bar{\rho}\bar{u} &0   \\
 \bar{v}^2+\frac{\gamma-1 }{2\gamma}(\bar{E}-2\bar{u}-\bar{v}^2)  & \frac{1-\gamma }{\gamma }\bar{\rho} & \frac{\gamma+1}{\gamma}\bar{\rho} \bar{v}   & \frac{\gamma -1}{2\gamma}\bar{\rho}   \\
  \bar{v} \bar{E} & 0  & \bar{\rho}\bar{E}   &\bar{\rho}\bar{v}
\end{pmatrix}.
\end{equation*}
Those are the solutions to the following equation:
\begin{equation}\label{0.4}
    (\lambda-\bar{v})(\lambda-\bar{v})[\lambda^2-2\lambda \bar{v}+\bar{v}^2-\bar{c}^2]=0,
\end{equation}
where $\bar{c}$ is the non-dimensional sonic speed given by 
\begin{equation*}
\bar{c}^2=\frac{\gamma\bar{p}}{\bar{\rho}}=\frac{\gamma-1}{2}(\bar{E}-2\bar{u}-\bar{v}^2).
\end{equation*}
By directly computation, (\ref{0.1}) admits four eigenvalues, i.e.,
\begin{equation}\label{0.5}
    \lambda_1=\bar{v}-\bar{c},\quad \lambda_2=\lambda_3=\bar{v},\quad  \lambda_4=\bar{v}+ \bar{c}.
\end{equation}
Thus, for $\gamma>1,~\bar{\rho}\ne 0,~\bar{p}\ne 0$, the system (\ref{0.1}) owns four real eigenvalues (containing a double eigenvalue); hence it is non-strictly hyperbolic. The corresponding eigenvectors are
\begin{equation}
\begin{aligned}
    \label{0.7}
    \mathbf{r}_1&=(-\frac{\bar\rho}{\bar{c}},\bar{c}-\bar{v},1,0)^{\top},\\
    \mathbf{r}_2&=(\frac{2\bar{\rho}}{\bar{E}-2\bar{u}-\bar{v}^2},1,0,0)^{\top},\\
   \mathbf{r}_3&=(-\frac{\bar{\rho}}{\bar{E}-2\bar{u}-\bar{v}^2},0,0,1)^{\top} ,\\
    \mathbf{r}_4&=(\frac{\bar\rho}{\bar{c}},-\bar{v}-\bar{c},1,0)^{\top},
    \end{aligned}
\end{equation}
which satisfy
\begin{equation}\label{0.8}
\nabla_{\bar{U}}\lambda_2(\bar{U})\cdot\mathbf{r}_2=\nabla_{\bar{U}}\lambda_3(\bar{U})\cdot\mathbf{r}_3\equiv 0,
\end{equation}
and
\begin{equation}\label{0.9}
    \begin{aligned}   \nabla_{\bar{U}}\lambda_1(\bar{U})\cdot\mathbf{r}_1&=(\gamma-1)\bar{c}+1\not\equiv 0,\\
      \nabla_{\bar{U}}\lambda_4(\bar{U})\cdot\mathbf{r}_4&=1+(\gamma-1)\bar{c}\not\equiv 0.
    \end{aligned}
\end{equation}
So we infer that the second and third characteristic fields are linearly degenerate, while the first and fourth characteristic fields are genuinely nonlinear.

\section{On the factor $\tau$ when comparing $\bar{\varrho}^{(\tau)}$ and $\bar{\varrho}$}\label{appendixfactor}

 This appendix aims to explain the reason why the factor $\tau$ appears in \eqref{2.61} and \eqref{3.70}. 
  
Since Radon measure solutions are defined in the sense of distribution, we demonstrate that the dimensionless transformation of solutions to the original boundary value problem in the measure framework introduces a factor of $\tau$ not present in the transformed differential equations (valid only for functions).

Without loss of generality, we consider Problem A presented in Section \ref{sec 2.1}, and rewrite the governing equations \eqref{2D-Euler eq}--\eqref{state eq} under the dimensionless transformation given by \eqref{2D scaling}--\eqref{2D scaling E}. Taking the mass conservation equation
\begin{equation}\label{mass eq}
\partial_x(\rho u) + \partial_y(\rho v) = 0
\end{equation}
as an example, we substitute the scaling \eqref{2D scaling}, yielding
\begin{equation*}
\partial_{\bar{x}}(\rho_{\infty}u_{\infty}\bar{\rho}(1 + \tau^2\bar{u})) + \frac{1}{\tau}\partial_{\bar{y}}(\rho_{\infty}u_{\infty}\bar{\rho}\tau\bar{v}) = 0.
\end{equation*}
After simplification, this reduces to
\begin{equation}\label{A.4}
\partial_{\bar{x}}(\bar{\rho}(1 + \tau^2\bar{u})) + \partial_{\bar{y}}(\bar{\rho}\bar{v}) = 0.
\end{equation}

We now consider instead the transformation of Eq. \eqref{mass eq} in the distributional sense. For any test function $\phi \in C_c^{1}(\overline{\Omega})$, it holds
\begin{align}\label{B.3}
\int_{\Omega}\phi \left[ \partial_x(\rho u) + \partial_y(\rho v) \right]\,\mathrm{d}x\mathrm{d}y = 0.
\end{align}
Using the equality
\begin{equation*}
\phi \left[ \partial_x(\rho u) + \partial_y(\rho v) \right] = \partial_x(\phi \rho u) + \partial_y(\phi \rho v) - \left( \rho u\partial_x\phi + \rho v\partial_y\phi \right)
\end{equation*}
and the Gauss-Green theorem, we derive
\begin{equation}\label{B.4}
\begin{aligned}
& \int_0^{\infty}\rho u\phi(x,\tau b(x))\mathrm{d}y + \int_{\infty}^{0}\rho_{\infty} u_{\infty}\phi(0, y)\mathrm{d}y - \int_0^{\infty}\rho v\phi(x,\tau b(x))\mathrm{d}x \\
& \quad - \int_{\Omega}\left( \rho u\partial_x\phi + \rho v\partial_y\phi \right)\mathrm{d}x\mathrm{d}y = 0.
\end{aligned}
\end{equation}
It follows from \eqref{2D scaling} that Eq. (\ref{B.4}) becomes
\begin{align*}
& \int_0^{\infty}\rho_{\infty}u_{\infty}\bar{\rho}(1 + \tau^2 \bar{u})\bar{\phi}^{\tau}(\bar{x}, b(\bar{x}))(\tau \mathrm{d}\bar{y}) + \int_{\infty}^{0}\rho_{\infty} u_{\infty}\bar{\phi}^{\tau}(0, \bar{y})(\tau\mathrm{d}\bar{y}) \\
& \quad - \int_0^{\infty}\tau\rho_{\infty}u_{\infty}\bar{\rho}\bar{v}\bar{\phi}^{\tau}(\bar{x}, b(\bar{x}))\mathrm{d}\bar{x} \\
& \qquad - \int_{\Omega'}\left( \rho_{\infty}u_{\infty}\bar{\rho}(1 + \tau^2 \bar{u})\partial_{\bar{x}}\bar{\phi}^{\tau} + \rho_{\infty}u_{\infty}\bar{\rho}\bar{v}\partial_{\bar{y}}\bar{\phi}^{\tau} \right)(\tau\mathrm{d}\bar{x}\mathrm{d}\bar{y}) = 0,
\end{align*}
where $\bar{\phi}^{(\tau)} \in C_c^{1}(\overline{\Omega'})$ denotes the test function in the $(\bar{x},\bar{y})$-plane. A simplification yields 
\begin{equation}\label{A.5}
\begin{aligned}
& \int_0^{\infty}\tau \bar{\rho}(1 + \tau^2 \bar{u})\bar{\phi}^{\tau}(\bar{x}, b(\bar{x}))\mathrm{d}\bar{y} + \int_{\infty}^{0} \tau \bar{\phi}^{\tau}(0, \bar{y})\mathrm{d}\bar{y} - \int_0^{\infty}\tau \bar{\rho}\bar{v}\bar{\phi}^{\tau}(\bar{x}, b(\bar{x}))\mathrm{d}\bar{x}  \\
& \quad - \int_{\Omega'}\tau\left( \bar{\rho}(1 + \tau^2 \bar{u})\partial_{\bar{x}}\bar{\phi}^{\tau} + \bar{\rho}\bar{v}\partial_{\bar{y}}\bar{\phi}^{\tau} \right)\mathrm{d}\bar{x}\mathrm{d}\bar{y} = 0. 
\end{aligned}
\end{equation}
It can be seen that all the terms in \eqref{A.5} contain a factor of $\tau$. Thus the dimensionless mass-conservation equation corresponding to  classical solutions takes the form
\begin{equation}\label{A.6}
\tau\left[ \partial _{\bar{x}}(\bar{\rho}( 1 + \tau^2\bar{u})) + \partial _{\bar{y}}(\bar{\rho} \bar{v}) \right] = 0.
\end{equation}
Compared with \eqref{A.4}, the two expressions differ by a factor of $\tau$. Since we consider the limit case as $\tau \to 0$, the Problem B, as a re-scaled problem, corresponds to an amplified treatment of Problem A. This explains why the density term vanishes at the boundary of the original problem as $\tau \to 0$ (corresponding to the case without a physical boundary), while in the small-disturbance equations, we observe a concentrated measure $\bar{\varrho}$ at the boundary.

Similarly, by the chain-rule, after dimensionless transformation of the momentum and energy equations, Eqs.  \eqref{2D-Euler eq} become
\begin{equation}\label{A.7}
\left\{
\begin{array}{l} 
\partial _{\bar{x}}(\bar{\rho}( 1 + \tau^2\bar{u}) )+ \partial _{\bar{y}}(\bar{\rho} \bar{v}) = 0, \\  
\tau^2 \left[ \partial _{\bar{x}}(\bar{\rho}\bar{u}( 1 + \tau^2\bar{u}) + \bar{p}) + \partial _{\bar{y}}(\bar{\rho} \bar{u}\bar{v}) \right] = 0, \\
\tau \left[ \partial _{\bar{x}}(\bar{\rho}\bar{v}( 1 + \tau^2\bar{u})) + \partial _{\bar{y}}(\bar{\rho}\bar{v}^2 + \bar{p}) \right] = 0, \\
\tau^2 \left[ \partial _{\bar{x}}(\bar{\rho}(1 + \tau ^2\bar{u})\bar{E}) + \partial _{\bar{y}}(\bar{\rho}\bar{v} \bar{E}) \right] = 0.
\end{array}
\right. 
\end{equation}
However, using the change-of-variable formula for multiple integrals, after scaling \eqref{2D scaling}--\eqref{2D scaling E}, \eqref{A.7} can be rewritten as
\begin{equation}\label{A.8}
\left\{
\begin{array}{l} 
\tau \left[ \partial _{\bar{x}}(\bar{\rho}( 1 + \tau^2\bar{u})) + \partial _{\bar{y}}(\bar{\rho} \bar{v})\right] = 0 , \\  
\tau^3 \left[ \partial _{\bar{x}}(\bar{\rho}\bar{u}( 1 + \tau^2\bar{u}) + \bar{p}) + \partial _{\bar{y}}(\bar{\rho} \bar{u}\bar{v}) \right] = 0, \\
\tau^2 \left[ \partial _{\bar{x}}(\bar{\rho}\bar{v}( 1 + \tau^2\bar{u})) + \partial _{\bar{y}}(\bar{\rho}\bar{v}^2 + \bar{p}) \right] = 0, \\
\tau^3 \left[ \partial _{\bar{x}}(\bar{\rho}(1 + \tau ^2\bar{u})\bar{E}) + \partial _{\bar{y}}(\bar{\rho}\bar{v} \bar{E}) \right] = 0
\end{array}
\right. 
\end{equation}
in the sense of distribution, where the two systems differ by a factor of  $\tau$. The same approach can be extended to the three-dimensional problem discussed in Section \ref{sec 4}.

\section*{Acknowledgments}
This work is supported by the Science and Technology Commission of Shanghai Municipality under Grants No.\,24ZR1420000 and No.\,22DZ2229014; Natural Science Foundation of Hubei under Grant No.\,2024AFB007.

%\bibliographystyle{plain} % for numbered citation & references
%\bibliography{ref}

\end{document}